\documentclass[12pt]{article}

\usepackage{color}
\usepackage{amssymb}
\usepackage{amsmath}
\usepackage{cite}

\begin{document}

\renewcommand{\citeleft}{{\rm [}}
\renewcommand{\citeright}{{\rm ]}}
\renewcommand{\citepunct}{{\rm,\ }}
\renewcommand{\citemid}{{\rm,\ }}

\newcounter{abschnitt}
\newtheorem{satz}{Theorem}
\newtheorem{theorem}{Theorem}[abschnitt]
\newtheorem{koro}[satz]{Corollary}
\newtheorem{coro}[theorem]{Corollary}
\newtheorem{prop}[theorem]{Proposition}
\newtheorem{lem}[theorem]{Lemma}
\newtheorem{conj}[theorem]{Conjecture}
\newtheorem{expls}[theorem]{Examples}

\renewenvironment{quote}{\list{}{\leftmargin=0.55in\rightmargin=0.55in}\item[]}{\endlist}

\newcounter{saveeqn}
\newcommand{\alpheqn}{\setcounter{saveeqn}{\value{abschnitt}}
\renewcommand{\theequation}{\mbox{\arabic{saveeqn}.\arabic{equation}}}}
\newcommand{\reseteqn}{\setcounter{equation}{0}
\renewcommand{\theequation}{\arabic{equation}}}

\hyphenpenalty=9000

\sloppy

\phantom{a}

\vspace{-1.7cm}

\begin{center}
\begin{Large} {\bf Projection Functions, Area Measures\\ and the Alesker-Fourier Transform} \\[0.6cm] \end{Large}

\begin{large} Felix Dorrek and Franz E. Schuster \end{large}
\end{center}

\vspace{-0.8cm}

\begin{quote}
\footnotesize{ \vskip 1cm \noindent {\bf Abstract.}
Dual to Koldobsky's notion of $j$-intersection bodies, the class of $j$-projection bodies is introduced, generalizing Minkowski's classical notion of projection bodies of convex bodies. A Fourier analytic characterization of $j$-projection bodies in terms of their area measures of order $j$ is obtained. In turn, this yields an equivalent characterization of $j$-projection bodies involving Alesker's Fourier type transform on translation invariant smooth spherical valuations. As applications of these results, several basic properties of $j$-projection bodies are established and new non-trivial examples are constructed. }
\end{quote}

\vspace{0.6cm}

\centerline{\large{\bf{ \setcounter{abschnitt}{1}
\arabic{abschnitt}. Introduction}}}

\alpheqn

\vspace{0.6cm}

The Busemann--Petty problem was one of the most famous problems in convex geometric analysis of the last century.
It asks whether the volume of an origin-symmetric convex body $K$ in $\mathbb{R}^n$ is smaller than that of another such body $L$, if all
central hyperplane sections of $K$ have smaller volume than those of $L$. (Here and throughout the article, it is assumed that $n \geq 3$.)
After more than 40 years and a long list of contributions it was
shown that the answer is affirmative if $n \leq 4$ and negative otherwise (see \textbf{\cite{gardner94, zhang99, gardkoldschlump}} and \linebreak the references therein).
The first crucial step in the final solution was taken by \linebreak Lutwak \textbf{\cite{lutwak84}} and later refined by Gardner \textbf{\cite{gardner94a}} who showed that the answer to the Busemann--Petty problem is affirmative if and only if every origin-symmetric convex body in $\mathbb{R}^n$ is an intersection body. This class of bodies first appeared in Busemann's definition of area in Minkowski geometry
and has attracted considerable attention in different subjects since the seminal paper by Lutwak (see, e.g., \textbf{\cite{Hab08, goodlutweil, GrinbZhang99, Koldobsky00, kaltonkold, Ludwig06}}
and the books \textbf{\cite{Gardner06, Koldobsky05, KoldobskyYaskin}}).

Since its final solution, several variants of the original Busemann--Petty problem have been investigated, each of which being related to a certain generalization
of the notion of intersection body in a similar way that Lutwak's intersection bodies are related to the Busemann--Petty problem (see \textbf{\cite{Koldobsky05, KoldobskyYaskin}}).
Of particular interest in this paper is the following notion of $j$-intersection bodies introduced by Koldobsky in 1999.

\vspace{0.3cm}

\noindent {\bf Definition} (\!\!\textbf{\cite{Koldobsky99}}) \emph{Let $1 \leq j \leq n - 1$ and let $D$ and $M$ be origin-symmetric star bodies in $\mathbb{R}^n$.
Then $D$ is called the \emph{$j$-intersection body} of $M$ if
\[\mathrm{vol}_{j}(D \cap E^{\bot}) = \mathrm{vol}_{n-j}(M \cap E)  \]
for every $n - j$ dimensional subspace $E$ of $\mathbb{R}^n$. The class of $j$-intersection bodies
is the closure in the radial metric of all $j$-intersection bodies of star bodies.}

\pagebreak

When $j=1$, the class of $1$-intersection bodies coincides with the closure of Lutwak's intersection bodies.
Also note that for $j > 1$, there may be star bodies $M$ for which a corresponding $j$-intersection body does not exist (cf.\ Theorem \ref{koldmain} below). However,
if $D$ is a $j$-intersection body of $M$ for some $1 \leq j \leq n - 1$, then $D$ is uniquely determined (see, e.g., \textbf{\cite[\textnormal{Corollary 7.2.7\,}]{Gardner06}}).

Since their definition by Koldobsky, $j$-intersection bodies have become objects of intensive investigations due to their connections to certain problems from functional analysis (\!\!\textbf{\cite{Koldobsky00, Schlieper07, Yaskin08}}), asymptotic geometric analysis (\!\!\textbf{\cite{KolPaoZym11}}), and complex geometry (\!\!\textbf{\cite{KolPaoZym13}}), as well as important variants of the Busemann--Petty problem
(\!\!\textbf{\cite{Koldobsky99, Milman06, Milman06, Yaskin14}}). The fundamental result on $j$-intersection bodies, which serves as starting point
for most subsequent investigations, is the following Fourier analytic characterization in terms of their radial functions.

\begin{satz} \label{koldmain} \emph{(\!\textbf{\cite{Koldobsky99, Koldobsky00}})}  Let $1 \leq j \leq n - 1$ and let $D$ and $M$ be origin-symmetric star bodies in $\mathbb{R}^n$. Then $D$ is the
$j$-intersection body of $M$ if and only if
\begin{equation*}
\mathbf{F}_{\!-j}\, \rho(D,\cdot)^j = \frac{(2\pi)^{n-j}j}{n-j}\,\rho(M,\cdot)^{n-j}.
\end{equation*}
\end{satz}

The operator $\mathbf{F}_{\!-j}$ denotes the (distributional) spherical Fourier transform of degree $-j$ originating in the work of Koldobsky (see Section 2 for details).
Note that the 'if' part of Theorem \ref{koldmain} is usually stated in the literature only for star bodies with smooth radial functions. However, the arguments
used in Section~4 of this paper show that this additional regularity assumption can be omitted.

\vspace{0.2cm}

Over the past decades, a remarkable correspondence between results about sections of star bodies through a fixed point and those concerning projections of convex bodies has asserted itself
(see, e.g., the books \textbf{\cite{Gardner06, Koldobsky05, Schneider14}}).
Thereby, the classical Brunn-Minkowski theory of convex bodies forms the ideal framework to deal with problems about projections while the \emph{dual} Brunn-Minkowski theory of star bodies provides the natural
setting for questions concerning sections. In this sense, the notion of $j$-intersection bodies and Theorem \ref{koldmain} belong to the latter dual theory. Surprisingly, so far no analogue of $j$-intersection bodies
has been thoroughly studied or even explicitly defined in the Brunn-Minkowski theory. In this article, we set out to remedy this neglect.

\vspace{0.3cm}

\noindent {\bf Definition} \emph{Let $1 \leq j \leq n - 1$ and let $K$ and $L$ be origin-symmetric convex bodies with non-empty interior in $\mathbb{R}^n$.
Then $K$ is called the \emph{$j$-projection body} of $L$ if
\[\mathrm{vol}_{j}(K|E^{\bot}) = \mathrm{vol}_{n-j}(L|E)  \]
for every $n - j$ dimensional subspace $E$ of $\mathbb{R}^n$. The class of $j$-projection bodies
is the closure in the Hausdorff metric of all $j$-projection bodies of convex bodies.}

\vspace{0.4cm}

When $j = 1$, the class of $1$-projection bodies coincides with the closure of Minkowski's projection bodies of convex bodies
which form a central notion in convex geometric analysis (see, e.g., \textbf{\cite{AbardiaBernig11, GrinbZhang99, Ludwig02, SchuWann12}} and the books \textbf{\cite{Gardner06, Schneider14}}).
We will see in Section 5 that, as in the case of $j$-intersection bodies, for $j > 1$, there exist convex bodies $L$ for which a corresponding
$j$-projection body does not exist. However, if $K$ is a $j$-projection body of $L$ for some $1 \leq j \leq n - 1$, then $K$ is uniquely determined (see, e.g., \textbf{\cite[\textnormal{Theorem 3.3.6}]{Gardner06}}).

Although for $j > 1$, the definition of $j$-projection bodies has not appeared before, special cases and examples
have been previously considered by several authors (see \textbf{\cite{McMullen84, McMullen87, Schneider96, Schnell94}}) and we recall them in Section 5. The main goal of this article, however,
is to start a systematic investigation of $j$-projection bodies of convex bodies. To this end, we not only establish a number of their basic properties, such as, invariance under non-degenerate linear transformations, but also obtain an array of new examples. These are based on our first main result which is the following Fourier analytic characterization.

\vspace{0.1cm}

\begin{satz} \label{mainthm1} Let $1 \leq j \leq n - 1$ and let $K$ and $L$ be origin-symmetric convex bodies with non-empty interior in $\mathbb{R}^n$. Then $K$ is the
$j$-projection body of $L$ if and only if
\begin{equation*}
\mathbf{F}_{\!-j}\,S_j(K,\cdot) = \frac{(2\pi)^{n-j}j}{(n-j)}\, S_{n-j}(L,\cdot).
\end{equation*}
\end{satz}

\vspace{0.1cm}

The Borel measures $S_j(K,\cdot)$, $1 \leq j \leq n - 1$, on $\mathbb{S}^{n-1}$ are Aleksandrov's \emph{area measures}
of the convex body $K$ (see Section 3 for details). Considering the still not fully understood correspondence
between results about sections and projections, we want to emphasize the astounding analogy between Theorem \ref{koldmain} and Theorem \ref{mainthm1},
where in order to pass from the characterization of $j$-intersection bodies to that of $j$-projection bodies,
certain powers of radial functions simply have to be replaced  by their 'dual' notion of area measures of the respective \linebreak orders.
We also note here that the case $j = 1$ of Theorem \ref{mainthm1} is equivalent to a previously known relation between
projection functions and the Fourier transform of surface area measures (see, e.g., \textbf{\cite{KolRyaZva04}}).

\vspace{0.2cm}

Theorem \ref{mainthm1} and a new relation between the spherical Fourier transform
and the Alesker-Fourier transform on spherical valuations (independently observed very recently by Goodey, Hug, and Weil \textbf{\cite{GooHugWei}})
lead to another characterization of $j$-projection bodies. In general, a \emph{valuation} on the space $\mathcal{K}^n$ of convex bodies in $\mathbb{R}^n$ is a map
$\phi: \mathcal{K}^n \rightarrow \mathcal{A}$ with values in an Abelian semigroup $\mathcal{A}$ such that
\[\phi(K) + \phi(L) = \phi(K \cup L) + \phi(K \cap L)   \]
whenever $K \cup L$ is convex. Since scalar valuations (where $\mathcal{A} = \mathbb{R}$ or $\mathbb{C}$) are on one hand intimately tied to the dissection theory of polytopes and on the other hand a generalization of measures, they have played an important role in discrete and integral geometry since the early 1900s (see, e.g., \textbf{\cite[\textnormal{Chapter 6}]{Schneider14}}).

\pagebreak

Recent work of Alesker \textbf{\cite{Alesker03}} and Bernig and Fu \textbf{\cite{BernigFu06}} uncovered natural product and convolution structures on the space $\mathbf{Val}^{\infty}$ of \emph{smooth} translation invariant valuations (precise definitions to follow) which were shown to encode important kinematic formulas of sectional and additive type, respectively. This not only prompted tremendous progress in integral geometry (see, e.g., \textbf{\cite{Bernig12,BernigFu10,BerFuSol15, Fu14}}) \linebreak but also aided in the resolution of several old mysteries. For example,
the Fourier type transform discovered by Alesker \textbf{\cite{Alesker11}},
\[\mathbb{F}: \mathbf{Val}^{\infty}_j \rightarrow \mathbf{Val}^{\infty}_{n-j}, \qquad 0 \leq j \leq n, \]
between spaces of smooth valuations of complementary degree intertwines
the product and convolution structures and, therefore, helped to explain the formal similarities between sectional and additive kinematic formulas. Moreover,
Bernig and Hug \textbf{\cite{BerHug15+}} made critical use of the Alesker-Fourier transform on \emph{spherical} valuations to obtain new kinematic formulas
for tensor valuations.

Spherical valuations correspond to \emph{spherical representations} of the group $\mathrm{SO}(n)$ (see Section 2 for details) and are not only crucial in the determination of integral geometric formulas for tensor valuations but also in the study of Minkowski valuations, that is, valuations with values in $\mathcal{A} = \mathcal{K}^n$ with Minkowski addition (see \textbf{\cite{SchuWann15, SchuWann16}}).
Let $\mathbf{Val}_j^{\infty,\mathrm{sph}}$ denote the subspace of smooth translation invariant spherical valuations of degree $j$.
Our second main result relates the Alesker-Fourier transform on even spherical valuations with $j$-projection bodies.

\begin{satz} \label{mainthm2} Let $1 \leq j \leq n - 1$ and let $K$ and $L$ be origin-symmetric convex bodies with non-empty interior in $\mathbb{R}^n$. Then $K$ is the
$j$-projection body of $L$ if and only if
\[\phi(K) = (\mathbb{F}\phi)(L)  \]
for all even $\phi \in \mathbf{Val}_j^{\infty,\mathrm{sph}}$.
\end{satz}

Note that Theorem \ref{mainthm2} is much easier to prove when the subspace $\mathbf{Val}_j^{\infty,\mathrm{sph}}$ is replaced by the entire space $\mathbf{Val}_j^{\infty}$.
The main point of Theorem \ref{mainthm2}, which also relates it to Theorem \ref{mainthm1}, is to consider spherical valuations only.

\vspace{0.2cm}

Instead of spherical scalar valuations let us now consider translation invariant and $\mathrm{SO}(n)$ equivariant even Minkowski valuations. Using Theorem \ref{mainthm1} or \ref{mainthm2}, it turns out that we can give a characterization of $j$-projection bodies in terms of \emph{a single pair} of such Minkowski valuations which are injective on origin-symmetric convex bodies and related by the Alesker-Fourier transform in a sense that we will make precise in Section 3. In the following corollary we exhibit one such pair explicitly, namely Minkowski's projection body operator of order $j$, $\Pi_j: \mathcal{K}^n \rightarrow \mathcal{K}^n$, and the (renormalized) mean section operator of Goodey and Weil \textbf{\cite{GooWei92, GooWei12, GooWei14}}, $\overline{\mathrm{M}}_{n-j}: \mathcal{K}^n \rightarrow \mathcal{K}^n$ (cf.\ Section 3 for definitions).

\begin{koro} \label{korotothm3} Let $1 \leq j \leq n - 1$ and let $K$ and $L$ be origin-symmetric convex bodies with non-empty interior in $\mathbb{R}^n$. Then $K$ is the
$j$-projection body of $L$ if and only if
\[\Pi_j K = \overline{\mathrm{M}}_{n-j}L.  \]
\end{koro}

The proofs of our main results will be presented in Section 4. The required background material from harmonic analysis, convex geometry, and the theory of valuations is the content of Sections 2 and 3.
In Section 5, we establish general properties of $j$-projection bodies, such as invariance under non-degenerate linear transformations or the fact that a polytope can only be the $j$-projection body of another polytope. We also review previously known examples in Section 5 and construct a large family of new ones.
In the final section, we relate $j$-intersection bodies and $j$-projection bodies via a duality transform motivated by Theorems \ref{koldmain} and \ref{mainthm1} and a celebrated result of Guan and Ma \textbf{\cite{GuaMa03}} on the Christoffel-Minkowski problem. We also discuss two open problems about $j$-projection bodies which are dual to recently resolved questions about $j$-intersection bodies.

\vspace{1cm}

\centerline{\large{\bf{ \setcounter{abschnitt}{2}
\arabic{abschnitt}. Background Material from Harmonic Analysis}}}

\reseteqn \alpheqn \setcounter{theorem}{0}

\vspace{0.6cm}

In the following we recall some well known facts about representations of the compact Lie group $\mathrm{SO}(n)$ and their
applications in the study of integral transforms of functions and measures on Grassmannians and the sphere.

Since the Lie group $\mathrm{SO}(n)$ is compact, all its irreducible representations are finite-dimensional. Moreover, the equivalence
classes of irreducible complex representations of $\mathrm{SO}(n)$ are uniquely determined by their highest weights (see, e.g., \textbf{\cite{knapp}}) which, in turn, can be
indexed by $\lfloor n/2 \rfloor$-tuples of integers $(\lambda_1,\lambda_2,\ldots,\lambda_{\lfloor n/2 \rfloor})$ such that
\begin{equation} \label{heiwei}
\left \{\begin{array}{ll} \lambda_1 \geq \lambda_2 \geq \ldots \geq \lambda_{\lfloor n/2 \rfloor} \geq 0 & \quad \mbox{for odd }n, \\
\lambda_1 \geq \lambda_2 \geq \ldots \geq \lambda_{n/2-1} \geq
|\lambda_{n/2}| & \quad \mbox{for even }n.
\end{array} \right .
\end{equation}

Throughout this article, we use $\bar{e} \in \mathbb{S}^{n-1}$ to denote an arbitrary but fixed point (the pole) of $\mathbb{S}^{n-1}$ and we write $\mathrm{SO}(n-1)$ for the stabilizer
in $\mathrm{SO}(n)$ of $\bar{e}$. \linebreak An important notion for our purposes is that of spherical representations of $\mathrm{SO}(n)$ with respect to $\mathrm{SO}(n-1)$.

\vspace{0.3cm}

\noindent {\bf Definition} \emph{Let $H$ be a closed subgroup of $\mathrm{SO}(n)$. A representation of $\mathrm{SO}(n)$ on a vector space $V$ is called \emph{spherical} with respect to $H$ if there exists an $H$-invariant non-zero $v \in V$, that is, $\vartheta v = v$ for every $\vartheta \in H$.}

\vspace{0.3cm}

The main result about spherical representations of a compact Lie group $G$ concerns the left regular representation
of $G$ on the Hilbert space $L^2(G/H)$ of square-integrable functions on the homogeneous space $G/H$ (see, \textbf{\cite[\textnormal{p.\ 17}]{takeuchi}}).

\pagebreak

\noindent However, we only require and state here the special case of this general result, where $G = \mathrm{SO}(n)$ and $H = \mathrm{SO}(n-1)$ and, consequently,
the homogeneous space $G/H$ is diffeomorphic to the sphere $\mathbb{S}^{n-1}$.

\begin{theorem}  \label{thmspher}
Every subrepresentation of the left regular representation of $\mathrm{SO}(n)$ on $L^2(\mathbb{S}^{n-1})$ is spherical with respect to $\mathrm{SO}(n-1)$.
Moreover, if $V$ is an $\mathrm{SO}(n)$ irreducible representation which is spherical with respect to $\mathrm{SO}(n-1)$, then $V$ is isomorphic to a subrepresentation of $L^2(\mathbb{S}^{n-1})$ and $\mathrm{dim}\,V^{\mathrm{SO}(n-1)} = 1$.
\end{theorem}

Here and in the following, we denote by $V^G$ the subspace of $G$-invariant vectors of a representation $V$ of a group $G$.

\begin{expls} \label{exp1} \end{expls}

\vspace{-0.2cm}

\begin{enumerate}
\item[(a)] The decomposition of $L^2(\mathbb{S}^{n-1})$ into an \emph{orthogonal} sum of $\mathrm{SO}(n)$ irreducible subspaces is given by
\begin{equation*} \label{decompsn1}
L^2(\mathbb{S}^{n-1}) = \bigoplus_{k \in \mathbb{N}} \mathcal{H}_k^n.
\end{equation*}
Here, $\mathcal{H}_k^n$ is the space of spherical harmonics of dimension $n$ and degree $k$. \linebreak
The highest weights associated with the spaces $\mathcal{H}_k^n$ are the $\lfloor n/2 \rfloor$-tuples
$(k,0,\ldots,0)$, $k \in \mathbb{N}$, and, by Theorem \ref{thmspher}, every irreducible representation of $\mathrm{SO}(n)$ which is spherical with respect to $\mathrm{SO}(n - 1)$ is isomorphic to one of the spaces $\mathcal{H}_k^n$.

A function or measure on $\mathbb{S}^{n-1}$ which is $\mathrm{SO}(n - 1)$ invariant is called \emph{zonal}. By Theorem \ref{thmspher},
each space $\mathcal{H}_k^n$ contains a $1$-dimensional subspace of zonal functions. This subspace
is spanned by the function $u \mapsto P_k^n(u \cdot \bar{e})$, where $P_k^n \in C([-1,1])$ is the \emph{Legendre polynomial} of dimension $n$ and degree $k$.
Letting $\pi_k: L^2(\mathbb{S}^{n-1}) \rightarrow \mathcal{H}_k^n$ denote the orthogonal projection, we write
\begin{equation} \label{fourierf}
f \sim \sum_{k=0}^{\infty} \pi_k f
\end{equation}
for the Fourier expansion of $f \in L^2(\mathbb{S}^{n-1})$. Recall that the Fourier series in (\ref{fourierf}) converges to $f$ in the $L^2$ norm and that
\begin{equation} \label{pikfLegPol}
(\pi_kf)(v) = N(n,k)\int_{\mathbb{S}^{n-1}}f(u)P_k^n(u\cdot v)\,du
\end{equation}
where $N(n,k)= \dim \mathcal{H}_k^n$ and integration is with respect to the $\mathrm{SO}(n)$ invariant probability measure on $\mathbb{S}^{n-1}$.

Since the orthogonal projection $\pi_k: L^2(\mathbb{S}^{n-1}) \rightarrow \mathcal{H}_k^n$ is self adjoint, it is consistent, by (\ref{pikfLegPol}), to extend it to the space $\mathcal{M}(\mathbb{S}^{n-1})$ of signed finite Borel measures by
\[(\pi_k \mu)(v) = N(n,k)\int_{\mathbb{S}^{n-1}}P_k^n(u\cdot v)\,d\mu(u).  \]

\pagebreak

It is not difficult to show that indeed $\pi_k\mu \in \mathcal{H}_k^n$ and that the Fourier expansion
\[\mu \sim \sum_{k=0}^{\infty} \pi_k \mu   \]
uniquely determines the measure $\mu \in \mathcal{M}(\mathbb{S}^{n-1})$.

\item[(b)] For $1\! \leq\! j\! \leq\! n - 1$, let $\mathrm{Gr}_{j,n}$ denote the Grassmann manifold of $j$-dimensional \linebreak subspaces of $\mathbb{R}^n$ and recall that
\begin{equation*} \label{grassmhomsp}
\mathrm{Gr}_{j,n} \cong \mathrm{SO}(n)/\mathrm{S}(\mathrm{O}(j) \times \mathrm{O}(n - j)).
\end{equation*}
The space $L^2(\mathrm{Gr}_{j,n})$ is a sum of orthogonal $\mathrm{SO}(n)$ irreducible subspaces with corresponding highest
weights $(\lambda_1,\ldots,\lambda_{\lfloor n/2 \rfloor})$ satisfying the following two conditions (see, e.g., \textbf{\cite[\textnormal{Theorem 8.49}]{knapp}}):
\begin{equation} \label{decompgri}
\left \{ \begin{array}{l} \lambda_k = 0\, \mbox{ for all } k > \min\{j, n-j\}, \phantom{wwwwwww} \\
\lambda_1, \ldots, \lambda_{\lfloor n/2 \rfloor}\, \mbox{ are all even.} \end{array} \right .
\end{equation}
Of particular importance for us is the subspace $L^2(\mathrm{Gr}_{j,n})^{\mathrm{sph}}$ of \emph{spherical functions} defined
as the orthogonal sum of all $\mathrm{SO}(n)$ irreducible subspaces in $L^2(\mathrm{Gr}_{j,n})$ which are spherical with respect to $\mathrm{SO}(n - 1)$.
By Theorem~\ref{thmspher}, Example \ref{exp1} (a), and (\ref{decompgri}),
\[L^2(\mathrm{Gr}_{j,n})^{\mathrm{sph}} = \bigoplus_{k \in \mathbb{N}} \Gamma_{(2k,0,\ldots,0)},   \]
where $\Gamma_{\lambda}$ denotes the $\mathrm{SO}(n)$ irreducible subspace of $L^2(\mathrm{Gr}_{j,n})$ of highest weight
$\lambda = (\lambda_1,\ldots,\lambda_{\lfloor n/2 \rfloor})$. Note that $L^2(\mathrm{Gr}_{j,n})^{\mathrm{sph}}$ is isomorphic as
$\mathrm{SO}(n)$ representation to the subspace of \emph{even} functions in $L^2(\mathbb{S}^{n-1})$.
\end{enumerate}

We now turn to convolution transforms of functions and measures on $\mathrm{SO}(n)$ and the homogeneous spaces $\mathbb{S}^{n-1}$ and $\mathrm{Gr}_{j,n}$.
These not only have the basic integral transforms we require, such as cosine and Radon transforms, as special cases but are also
crucial in our discussion of Minkowski valuations in the next section. In order to keep the exposition brief,
we identify integrable functions on $\mathrm{SO}(n)$, $\mathbb{S}^{n-1}$, or $\mathrm{Gr}_{j,n}$ with absolutely continuous measures with respect to the
corresponding $\mathrm{SO}(n)$ invariant probability measures and, thus, state most formulas for measures only.

First recall that the convolution $\mu \ast \sigma$ of signed measures $\mu, \sigma$ on $\mathrm{SO}(n)$ is defined as
the pushforward of the product measure $\mu \otimes \sigma$ by the group multiplication $m: \mathrm{SO}(n) \times \mathrm{SO}(n) \rightarrow \mathrm{SO}(n)$,
that is, $\mu \ast \sigma = m_*(\mu \otimes \sigma)$ or, more explicitly,
\[\int_{\mathrm{SO}(n)}\!\!\! f(\vartheta)\, d(\mu \ast \sigma)(\vartheta)=\int_{\mathrm{SO}(n)}\! \int_{\mathrm{SO}(n)}\!\!\! f(\eta \theta)\,d\mu(\eta)\,d\sigma(\theta), \qquad f \in C(\mathrm{SO}(n)).   \]

For a measure $\mu$ on $\mathrm{SO}(n)$, let $l_{\vartheta}\mu$ and $r_{\vartheta}\mu$ denote the pushforward of $\mu$ by
the left and right translations by $\vartheta \in \mathrm{SO}(n)$, respectively. We also often use $\vartheta\mu := l_{\vartheta}\mu$ for the left translation of $\mu$.
It follows from the definition of $\mu \ast \sigma$ that
\begin{equation} \label{roteq}
(l_{\vartheta} \mu) \ast \sigma = l_{\vartheta}(\mu \ast \sigma) \qquad \mbox{and} \qquad \mu \ast (r_{\vartheta} \sigma) = r_{\vartheta}(\mu \ast \sigma)
\end{equation}
for every $\vartheta \in G$. Moreover, the convolution of measures on $\mathrm{SO}(n)$ is associative but in general not commutative. In fact,
if $\mu, \sigma$ are measures on $\mathrm{SO}(n)$, then
\begin{equation*} \label{commute}
\widehat{\mu \ast \sigma}=\widehat{\sigma} \ast \widehat{\mu},
\end{equation*}
where $\widehat{\mu}$ denotes the pushforward of $\mu$ by the group inversion, that is,
\begin{equation*}
\int_{\mathrm{SO}(n)}\!\!\! f(\vartheta)\,d\widehat{\mu}(\vartheta) = \int_{\mathrm{SO}(n)}\!\!\! f(\vartheta^{-1})\,d\mu(\vartheta), \qquad f \in C(\mathrm{SO}(n)).
\end{equation*}

In order to define the convolution of measures on $\mathbb{S}^{n-1}$ and $\mathrm{Gr}_{j,n}$, we make use of the diffeomorphisms
\[\mathbb{S}^{n-1}=\mathrm{SO}(n)/\mathrm{SO}(n-1) \quad \mbox{and} \quad \mathrm{Gr}_{j,n}=\mathrm{SO}(n)/\mathrm{S}(\mathrm{O}(j) \times \mathrm{O}(n-j)).  \]
Indeed, if $H$ is a closed subgroup of $\mathrm{SO}(n)$, then there is a natural one-to-one \linebreak correspondence
between measures on $\mathrm{SO}(n)/H$ and right $H$-invariant measures on $\mathrm{SO}(n)$ (see, e.g., \textbf{\cite{GrinbZhang99, Schuster10}} for a detailed description).
Using this identification, the convolution of measures on $\mathrm{SO}(n)$ induces a convolution product of measures on $\mathrm{SO}(n)/H$ as follows:
Let $\pi: \mathrm{SO}(n) \rightarrow \mathrm{SO}(n)/H$ denote the canonical projection. The convolution of measures $\mu$ and $\sigma$ on $\mathrm{SO}(n)/H$ is defined by
\begin{equation} \label{defconvhomsp}
\mu \ast \sigma = \pi_* m_*(\pi^*\mu \otimes \pi^*\sigma),
\end{equation}
where $\pi_*$ and $\pi^*$ denote the pushforward and pullback by $\pi$, respectively. Note that, by (\ref{roteq}), definition (\ref{defconvhomsp}) is consistent with the identification
of measures on $\mathrm{SO}(n)/H$ with right $H$-invariant measures on $\mathrm{SO}(n)$. In the same way, the convolution of measures on different
homogeneous spaces can be defined: Let $H_1, H_2$ be two closed subgroups of $\mathrm{SO}(n)$ and denote by $\pi_i: \mathrm{SO}(n) \rightarrow \mathrm{SO}(n)/H_i$, $i = 1,2$, the respective projections.
If, say, $\mu$ is a measure on $\mathrm{SO}(n)/H_1$ and $\sigma$ a measure on $\mathrm{SO}(n)/H_2$, then
\[\mu \ast \sigma = \pi_{2*} m_*(\pi_1^*\mu \otimes \pi_2^*\sigma),   \]
defines a measure on $\mathrm{SO}(n)/H_2$.

Since the projection $\pi: \mathrm{SO}(n) \rightarrow \mathrm{SO}(n)/H$ is given by $\pi(\vartheta) = \vartheta \bar{E}$, where $H$ is the stabilizer in $G$ of $\bar{E} \in \mathrm{SO}(n)/H$ (note that we write $\bar{e}$ instead of $\bar{E}$ when $H = \mathrm{SO}(n - 1)$),
the convolution of a measure $\mu$ on $\mathrm{SO}(n)$ with the Dirac measure $\delta_{\bar{E}}$ on $\mathrm{SO}(n)/H$ yields
\begin{equation} \label{dirac}
\mu \ast \delta_{\bar{E}} = \int_H r_{\vartheta}\mu\,d\vartheta \qquad \mbox{and} \qquad \delta_{\bar{E}} \ast \mu  =\int_{H}l_{\vartheta} \mu\,d\vartheta.
\end{equation}

\pagebreak

Thus, $\delta_{\bar{E}} \ast \mu$ is left $H$-invariant, $\mu \ast \delta_{\bar{E}}$ is right $H$-invariant, and $\delta_{\bar{E}}$ is the unique rightneutral
element for the convolution of measures on $\mathrm{SO}(n)/H$. Generalizing the notion of zonal measures on $\mathbb{S}^{n-1}$, a left $H$-invariant measure on $\mathrm{SO}(n)/H$
is called zonal. If $\mu$ and $\sigma$ are measures on $\mathrm{SO}(n)/H$, then, by (\ref{dirac}),
\[\mu \ast \sigma = (\mu \ast \delta_{\bar{E}}) \ast \sigma = \mu \ast (\delta_{\bar{E}} \ast \sigma).  \]
Consequently, for the convolution of measures on $\mathrm{SO}(n)/H$, the right hand side measure can always assumed to be zonal.

Before we discuss important specific examples, we recall one more critical property of the convolution of measures on $\mathbb{S}^{n-1}$.
Using the identification of a zonal measure $\mu$ on $\mathbb{S}^{n-1}$ with a measure on $[-1,1]$ and the Funk-Hecke Theorem, one can show (cf.\ \textbf{\cite{Schuster07}})
that the Fourier expansion of $\sigma \ast \mu$ is given by
\begin{equation} \label{funkheckgen}
\sigma \ast \mu \sim \sum_{k=0}^{\infty} a_k^n[\mu]\,\pi_k\sigma,
\end{equation}
where the numbers
\begin{equation*}
a_k^n[\mu] = \omega_{n-1} \int_{-1}^1 P_k^n(t)\,(1-t^2)^{\frac{n-3}{2}}\,d\mu(t)
\end{equation*}
are called the \emph{multipliers} of the convolution transform $\sigma \mapsto \sigma \ast \mu$.
Here, $\omega_{n-1}$ is the surface area of the $(n-1)$-dimensional Euclidean unit ball.

\begin{expls} \label{exptransforms} \end{expls}

\vspace{-0.2cm}

\begin{enumerate}
\item[(a)] Let $1 \leq j \leq n - 1$ and let $|\cos(E,F)|$ denote the cosine of the angle between two subspaces $E, F \in \mathrm{Gr}_{j,n}$ (see, e.g., \textbf{\cite{GooZha98}}).
The {\it cosine transform} $\mathrm{C}_j\mu$ of a measure $\mu$ on $\mathrm{Gr}_{j,n}$ is the continuous function on $\mathrm{Gr}_{j,n}$ defined by
\[(\mathrm{C}_j\mu)(F) = \int_{\mathrm{Gr}_{j,n}} |\cos(E,F)|\,d\mu(E).  \]
It is not difficult to show that
\begin{equation} \label{cosconv}
\mathrm{C}_j \mu = \mu \ast |\cos(\bar{E},\,\cdot\,)|,
\end{equation}
where $\bar{E} \in \mathrm{Gr}_{j,n}$ again denotes the image of the identity under the projection $\pi: \mathrm{SO}(n) \rightarrow \mathrm{Gr}_{j,n}$.
In particular, the cosine transform is a linear and self-adjoint operator which is $\mathrm{SO}(n)$ equivariant and
maps smooth functions to smooth ones, that is,
\[\mathrm{C}_j: C^{\infty}(\mathrm{Gr}_{j,n}) \rightarrow C^{\infty}(\mathrm{Gr}_{j,n}).   \]
Moreover, since $|\cos(E,F)| = |\cos(E^{\bot},F^{\bot})|$, we have
\begin{equation} \label{cifbot}
(\mathrm{C}_j\mu)^{\bot} = \mathrm{C}_{n-j}\mu^{\bot}
\end{equation}
where $\mu^{\bot}:= \bot_*\mu$ denotes the pushforward of $\mu$ by the orthogonal complement map $\bot: \mathrm{Gr}_{j,n} \rightarrow \mathrm{Gr}_{n-j,n}$.

\pagebreak

It is a classical fact (see, e.g., \textbf{\cite[\textnormal{Chapter 3}]{Groemer}}) that the {\it spherical} cosine transform
$\mathrm{C}_1$ is injective and, thus, by (\ref{cifbot}), so is $\mathrm{C}_{n-1}$.
For $1 < j < n - 1$, Goodey and Howard \textbf{\cite{goodeyhoward}} first showed that the cosine transform $\mathrm{C}_j$ is \emph{not}
injective. A precise description of its kernel was given by Alesker and Bernstein \textbf{\cite{AlBern}}.
However, Goodey and Zhang \textbf{\cite[\textnormal{Lemma 2.1}]{GooZha98}} proved that the restriction of $\mathrm{C}_j$
to spherical functions in $L^2(\mathrm{Gr}_{j,n})^{\mathrm{sph}}$ is injective and, moreover, when restricted to the subspace of smooth spherical functions
\begin{equation} \label{smoothspherical}
C^{\infty}(\mathrm{Gr}_{j,n})^{\mathrm{sph}} := \mathrm{cl}_{C^{\infty}} \bigoplus_{k \in \mathbb{N}} \Gamma_{(2k,0,\ldots,0)}
\end{equation}
the cosine transform $\mathrm{C}_j$ is {\it bijective} for every $1 \leq i \leq n - 1$. Here, $\mathrm{cl}_{C^{\infty}}$ denotes the
closure in the $C^{\infty}$ topology.

\item[(b)] Let $1 \leq i \neq j \leq n - 1$. For $F \in \mathrm{Gr}_{j,n}$, we write $\mathrm{Gr}_{i,n}^F$ for the submanifold of
$\mathrm{Gr}_{i,n}$ which comprises of all $E \in \mathrm{Gr}_{i,n}$ that contain (respectively, are contained in) $F$. The
{\it Radon transform} $\mathrm{R}_{i,j}: L^2(\mathrm{Gr}_{i,n}) \rightarrow L^2(\mathrm{Gr}_{j,n})$ is defined by
\[(\mathrm{R}_{i,j}f)(F) = \int_{\mathrm{Gr}_{i,n}^F}f(E)\,d\nu_i^F(E),   \]
where $\nu_i^F$ is the unique invariant probability measure on $\mathrm{Gr}_{i,n}^F$.
It is well known that, for $1 \leq i < j < k \leq n - 1$,  we have
\[\mathrm{R}_{i,k} = \mathrm{R}_{j,k} \circ \mathrm{R}_{i,j} \qquad \mbox{and} \qquad \mathrm{R}_{k,i} = \mathrm{R}_{j,i} \circ \mathrm{R}_{k,j}   \]
and that $\mathrm{R}_{j,i}$ is the \emph{adjoint} of $\mathrm{R}_{i,j}$. Using this latter fact, one can define the Radon transform of a measure $\mu$ on $\mathrm{Gr}_{i,n}$ by
\[\int_{\mathrm{Gr}_{j,n}}f(F)\,d(\mathrm{R}_{i,j}\mu)(F) = \int_{\mathrm{Gr}_{i,n}} (\mathrm{R}_{j,i}f)(E)\,d\mu(E), \qquad f \in C(\mathrm{Gr}_{j,n}).\]
Also the Radon transform intertwines the orthogonal complement map. More precisely,
\begin{equation} \label{rijbot}
(\mathrm{R}_{i,j}\mu)^{\bot} = \mathrm{R}_{n - i,n - j}\mu^{\bot}.
\end{equation}

For $1 \leq i < j \leq n - 1$ let $\lambda_{i,j}$ denote the probability measure on $\mathrm{Gr}_{j,n}$ which is uniformly concentrated on the submanifold
\[\{\vartheta \bar{E} \in \mathrm{Gr}_{j,n}: \vartheta \in \mathrm{S}(\mathrm{O}(i) \times \mathrm{O}(n - i))\}.   \]
It is not difficult to show (see, e.g., \textbf{\cite{GrinbZhang99}}) that for measures $\mu$ on $\mathrm{Gr}_{i,n}$ and $\nu$ on $\mathrm{Gr}_{j,n}$, we have
\begin{equation} \label{radconv}
\mathrm{R}_{i,j}\mu = \mu \ast \lambda_{i,j} \qquad \mbox{and} \qquad \mathrm{R}_{j,i}\nu = \nu \ast \widehat{\lambda}_{i,j}.
\end{equation}
In particular, the Radon transform is a linear $\mathrm{SO}(n)$ equivariant operator which maps smooth functions to smooth ones, that is,
\[\mathrm{R}_{i,j}: C^{\infty}(\mathrm{Gr}_{i,n}) \rightarrow C^{\infty}(\mathrm{Gr}_{j,n}).   \]

\pagebreak

It follows from results of Grinberg \textbf{\cite{grinberg86}} that if $1 \leq i < j \leq n - 1$, then $\mathrm{R}_{i,j}$ is injective if and
only if $i + j \leq n$, whereas if $i > j$, then $\mathrm{R}_{i,j}$ is injective if and only if $i + j \geq n$. Moreover, Goodey and Zhang
\textbf{\cite{GooZha98}} proved that for all $1 \leq i \neq j \leq n -1$ the restriction of the Radon transform $\mathrm{R}_{i,j}$ to spherical functions is injective
and that
\[\mathrm{R}_{i,j}: C^{\infty}(\mathrm{Gr}_{i,n})^{\mathrm{sph}} \rightarrow C^{\infty}(\mathrm{Gr}_{j,n})^{\mathrm{sph}}   \]
is a \emph{bijection}.

\item[(c)] Let $\mathcal{S}(\mathbb{R}^n)$ denote the Schwartz space of complex valued, rapidly decreasing, infinitely differentiable \emph{test functions} on $\mathbb{R}^n$
endowed with its standard topology (see, e.g., \textbf{\cite[\textnormal{Chapter 2.5}]{Koldobsky05}}). We call a linear, continuous functional on $\mathcal{S}(\mathbb{R}^n)$ a \emph{distribution} over $\mathcal{S}(\mathbb{R}^n)$. Note that any locally integrable function on $\mathbb{R}^n$ satisfying
a power growth condition at infinity (cf.\ \textbf{\cite[\textnormal{p.\ 34}]{Koldobsky05}}) determines a distribution acting by integration.

The \emph{Fourier transform} $\mathrm{F}: \mathcal{S}(\mathbb{R}^n) \rightarrow \mathcal{S}(\mathbb{R}^n)$ is defined by
\[(\mathrm{F}\tau)(x) = \int_{\mathbb{R}^n} \tau(y)\exp(-i\,x\cdot y)\,dy.  \]
It is well known that $\mathrm{F}$ is an $\mathrm{SO}(n)$ equivariant (topological) isomorphism of the Schwartz space $\mathcal{S}(\mathbb{R}^n)$.
Moreover, $\mathrm{F}$ is self-adjoint on $\mathcal{S}(\mathbb{R}^n)$. This motivates the definition of the Fourier transform $\mathrm{F}\nu$ of a
distribution $\nu$ over $\mathcal{S}(\mathbb{R}^n)$ as the distribution acting by
\[\langle \mathrm{F}\nu,\tau \rangle = \langle \nu, \mathrm{F}\tau\rangle, \qquad \tau \in \mathcal{S}(\mathbb{R}^n).  \]
A distribution $\nu$ over $\mathcal{S}(\mathbb{R}^n)$ is called \emph{even homogeneous of degree $p \in \mathbb{R}$} if
\[\langle \nu, \tau(\,\cdot\,/\lambda) \rangle = |\lambda|^{n+p} \langle \nu,\tau \rangle  \]
for every $\tau \in \mathcal{S}(\mathbb{R}^n)$ and every $\lambda \in \mathbb{R}\backslash \{0\}$.
For the rest of this article, we only consider even homogeneous distributions $\nu$. Note that in this case,
\begin{equation} \label{Fsquaredeven}
\mathrm{F}^2\nu = (2\pi)^n\nu.
\end{equation}

Moreover, Koldobsky \textbf{\cite[\textnormal{Lemma 2.21}]{Koldobsky05}} observed the following crucial fact.

\begin{lem} \label{fourhom} The Fourier transform of an even homogeneous distribution of degree $p$ is an even homogeneous distribution of degree $-n-p$.
\end{lem}

Now consider the space $C_e^{\infty}(\mathbb{S}^{n-1})$ of all real-valued \emph{even} smooth functions on $\mathbb{S}^{n-1}$ endowed with its standard
Fr\'echet space topology. For $p > -n$ and $f \in C_e^{\infty}(\mathbb{S}^{n-1})$, we denote by $f_p$ the homogeneous extension of $f$ of degree $p$ to
$\mathbb{R}^n\backslash\{0\}$, that is,
\[f_p(x) = \|x\|^pf\left (\frac{x}{\|x\|} \right ), \qquad x \in \mathbb{R}^n\backslash \{0\}.  \]

\pagebreak

Since $p > -n$, $f_p$ is locally integrable and determines an even homogeneous distribution of degree $p$ acting on test functions by integration.
Thus, by Lemma \ref{fourhom}, $\mathrm{F}f_p$ is an even homogeneous distribution of degree $-n - p$. It was first noted in \textbf{\cite{GooYasYas11}} that, for $-n < p < 0$,
$\mathrm{F}f_p$ is, in fact, an infinitely differentiable function on $\mathbb{R}^n \backslash\{0\}$ (which is even and homogeneous of degree $-n - p$).
This gives rise to an operator $\mathbf{F}_p$ on $C_e^{\infty}(\mathbb{S}^{n-1})$, called the \emph{spherical Fourier transform of degree $p \in (-n,0)$}, defined by
\[\mathbf{F}_pf = \mathrm{F}f_p|_{\mathbb{S}^{n-1}}, \qquad f \in C_e^{\infty}(\mathbb{S}^{n-1}).   \]
Clearly, $\mathbf{F}_p$ is a linear and $\mathrm{SO}(n)$ equivariant map. Hence, by Schur's lemma, $\mathbf{F}_p$ acts as a
multiplier transformation on the spaces $\mathcal{H}_{2k}^n$, $k \in \mathbb{N}$. Its multipliers $a_{2k}^n[\mathbf{F}_p]$ were
determined in \textbf{\cite{GooYasYas11}} and are given by
\begin{equation} \label{multsphfour}
a_{2k}^n[\mathbf{F}_p] = \pi^{n/2}2^{n+p}(-1)^k\frac{\Gamma\left (\frac{2k + n + p}{2}\right )}{\Gamma \left (\frac{2k - p}{2} \right )}.
\end{equation}
We will give a convolution representation of $\mathbf{F}_p$ at the end of the next section. For now, we just note that (\ref{multsphfour}) implies that for $p \in (-n,0)$,
\[\mathbf{F}_p: C_e^{\infty}(\mathbb{S}^{n-1}) \rightarrow C_e^{\infty}(\mathbb{S}^{n-1})  \]
is bijective and that, by Lemma \ref{fourhom}, (\ref{Fsquaredeven}), and the definition of $\mathbf{F}_p$,
\begin{equation} \label{sphfourinv}
\mathbf{F}_{-n-p}(\mathbf{F}_pf) = (2\pi)^n f
\end{equation}
holds for every $f \in C_e^{\infty}(\mathbb{S}^{n-1})$. Moreover, as a multiplier transformation $\mathbf{F}_p$ is self-adjoint and, hence,
admits an extension to the space $C^{-\infty}_e(\mathbb{S}^{n-1})$ of continuous, linear functionals on $C^{\infty}_e(\mathbb{S}^{n-1})$ defined by
\[\langle \mathbf{F}_p \nu, f \rangle = \langle \nu,  \mathbf{F}_p f \rangle  \]
for $\nu \in C^{-\infty}_e(\mathbb{S}^{n-1})$ and $f \in C_e^{\infty}(\mathbb{S}^{n-1})$. Since every even (signed) Borel measure $\mu$ on $\mathbb{S}^{n-1}$
determines an element $\nu_{\mu} \in C^{-\infty}_e(\mathbb{S}^{n-1})$ by integration, we use the continuous, linear injection $\mu \mapsto \nu_{\mu}$, to identify the space $\mathcal{M}_e(\mathbb{S}^{n-1})$ of even Borel measures with a subspace of $C^{-\infty}_e(\mathbb{S}^{n-1})$. In this way, the spherical Fourier transform $\mathbf{F}_{p}\mu$ of
$\mu \in \mathcal{M}_e(\mathbb{S}^{n-1})$ is defined.

Finally, we state a fundamental relation between the spherical Fourier transform and certain Radon transforms which was first observed by Koldobsky \textbf{\cite{Koldobsky99}}
(see also \textbf{\cite{Milman06}}). Here, $\kappa_m$ is the $m$-dimensional volume of the Euclidean unit ball in $\mathbb{R}^m$.

\begin{prop} \label{spherfourradon} Suppose that $1 \leq j \leq n - 1$. Then
\[\mathrm{R}_{1,n-j} \circ \mathbf{F}_{-j} = \frac{(2\pi)^{n-j}\,j\,\kappa_j}{(n-j)\,\kappa_{n-j}}\,\bot_* \circ \mathrm{R}_{1,j}.  \]
\end{prop}

Note that here and in the following, we identify $C^{\infty}_e(\mathbb{S}^{n-1})$ with $C^{\infty}(\mathrm{Gr}_{1,n})$ and
the transform $\mathrm{R}_{i,i}$, $1 \leq i \leq n - 1$, with the identity map.

\end{enumerate}

\pagebreak

\centerline{\large{\bf{ \setcounter{abschnitt}{3}
\arabic{abschnitt}. Area Measures and Valuation Theory}}}

\reseteqn \alpheqn \setcounter{theorem}{0}

\vspace{0.6cm}

In this section, we collect background material from convex geometry, mainly about area measures and the Christoffel-Minkowski problem.
We also recall the material from the theory of valuations required for the proofs of Theorem \ref{mainthm2} and Corollary \ref{korotothm3}
and, in particular, note a critical relation between the spherical Fourier transform and the Alesker-Fourier transform of spherical valuations.
For most of this background material we recommend the book by Schneider \textbf{\cite{Schneider14}}.

First recall that $\mathcal{K}^n$ denotes the space of convex bodies (that is, of non-empty, compact, convex sets) in $\mathbb{R}^n$
and that we always assume that $n \geq 3$. A convex body $K \in \mathcal{K}^n$ is uniquely determined by the values of its support function
$h(K,u)=\max \{u \cdot x: x \in K\}$ for $u \in \mathbb{S}^{n-1}$. Note that, for origin-symmetric $K$,
\begin{equation} \label{suppvol1}
\mathrm{vol}_1(K|\mathrm{span}\{u\}) = 2h(K,u).
\end{equation}

A compact set $L$ in $\mathbb{R}^n$ which is star-shaped with respect to the origin is uniquely determined by its
radial function $\rho(L,u)=\max\{\lambda \geq 0: \,\lambda u\in L\}$ for $u \in \mathbb{S}^{n-1}$. If $\rho(L,\cdot)$ is positive and
continuous, we call $L$ a star body. If a convex body $K \in \mathcal{K}^n$ contains the origin in its interior, then
\begin{equation} \label{suprad}
\rho(K^*,\cdot)=h(K,\cdot)^{-1} \qquad \mbox{and} \qquad h(K^*,\cdot)=\rho(K,\cdot)^{-1},
\end{equation}
where $K^* = \{x \in \mathbb{R}^n: x \cdot y \leq 1 \mbox{ for all } y \in K\}$ is the polar body of $K$.

A classical result of Minkowski states that the volume of a Minkowski linear combination $\lambda_1K_1 + \cdots + \lambda_mK_m$,
where $K_1, \ldots, K_m \in \mathcal{K}^n$ and $\lambda_1, \ldots, \lambda_m \geq 0$, can be expressed as a homogeneous polynomial of
degree $n$, that is,
\begin{equation*} \label{mixed}
V_n(\lambda_1K_1 + \cdots +\lambda_m K_m)=\sum \limits_{1\leq j_1,\ldots, j_n\leq m} V(K_{j_1},\ldots,K_{j_n})\lambda_{j_1}\cdots\lambda_{j_n}.
\end{equation*}
The symmetric coefficients $V(K_{j_1},\ldots,K_{j_n})$ are called the {\it mixed volumes} of $K_{j_1}, \ldots, K_{j_n}$. For $K, L \in \mathcal{K}^n$ and $0 \leq j \leq n$, we denote
the mixed volume with $j$ copies of $K$ and $n - j$ copies of $L$ by $V(K[j],L[n-j])$ and we write $V_j(K)$ for the $j$th \emph{intrinsic volume} of $K$ defined by
\[\kappa_{n-j}V_j(K) = {n \choose j} V(K[j],B[n-j]).  \]

The {\it surface area measure} $S_{n-1}(K,\cdot)$ of a convex body $K \in \mathcal{K}^n$ is the Borel measure on $\mathbb{S}^{n-1}$ defined, for $\omega \subseteq \mathbb{S}^{n-1}$, as
the $(n-1)$-dimensional Hausdorff measure of the set of all boundary points of $K$ at which there exists a normal vector of
$K$ belonging to $\omega$. If $B$ denotes the Euclidean unit ball in $\mathbb{R}^n$, then, for every $r \geq 0$, the surface area measure of $K$ satisfies the
Steiner-type formula
\begin{equation*} \label{steiner}
S_{n-1}(K + r B,\cdot)=\sum \limits_{j=0}^{n-1} r^{n-1-j} {n - 1 \choose j}S_j(K,\cdot).
\end{equation*}
The Borel measure $S_j(K,\cdot)$, $0 \leq j \leq n - 1$, on $\mathbb{S}^{n-1}$ is called the
{\it area measure of order $j$} of $K \in \mathcal{K}^n$. It is uniquely determined by the property that
\begin{equation} \label{defsi}
V(K[j],B[n-1-j],L) = \frac{1}{n}\int_{\mathbb{S}^{n-1}} h(L,u)\,dS_j(K,u)
\end{equation}
for all $L \in \mathcal{K}^n$. If $K \in \mathcal{K}^n$ has non-empty interior, then, by a theorem of
Aleksandrov-Fenchel-Jessen (see, e.g., \textbf{\cite[\textnormal{p.\ 449}]{Schneider14}}), each of the
measures $S_j(K,\cdot)$, $1 \leq j \leq n - 1$, determines $K$ up to translations. In particular,
if $K$ is origin-symmetric, then $S_j(K,\cdot)$ is an even measure on $\mathbb{S}^{n-1}$ and, thus, can be
identified with a measure on $\mathrm{Gr}_{1,n}$ of the same total mass. Using this identification, the important
\emph{Cauchy-Kubota formula} can be stated as follows: For every $1 \leq j \leq n - 1$ and origin-symmetric $K \in \mathcal{K}^n$,
\begin{equation} \label{radsurfcos}
(\bot_* \circ \mathrm{R}_{j,n-1})\mathrm{vol}_j(K|\,\cdot\,) = \frac{\kappa_j}{2\kappa_{n-1}}\, \mathrm{C}_1S_j(K,\cdot) = \frac{\kappa_j}{\kappa_{n-1}}\, V_j(K|\,\cdot\,)^{\bot}.
\end{equation}

In view of Proposition \ref{spherfourradon} and (\ref{radsurfcos}), the following result about the injectivity of Radon transforms of
projection functions is important for our purpose.

\begin{prop}[\!\!\cite{Aleksandrov37, Goodey98}] \label{goodeyaleks} Suppose that $1 \leq i, j \leq n - 1$ and let $K, L \in \mathcal{K}^n$ be origin-symmetric and have non-empty interior.
If $\mathrm{R}_{i,j}\mathrm{vol}_i(K|\,\cdot\,) = \mathrm{R}_{i,j}\mathrm{vol}_i(L|\,\cdot\,)$ on $\mathrm{Gr}_{j,n}$, then $K = L$.
\end{prop}

Recall that a body $K \in \mathcal{K}^n$ is called a \emph{zonoid} if it can be approximated in the Hausdorff metric by finite Minkowski sums of line segments.
Let $\mathcal{Z}^n_s$ denote the class of \emph{origin-symmetric zonoids} in $\mathbb{R}^n$. It is well known that a convex body $K$ belongs to the class $\mathcal{Z}^n_s$ if and only if
its support function can be represented in the form
\begin{equation} \label{zonoidchar}
h(K,\cdot) = \mathrm{C}_1 \mu_K
\end{equation}
for some uniquely determined non-negative $\mu_K \in \mathcal{M}_e(\mathbb{S}^{n-1})$.

For $1 \leq j \leq n - 1$, Weil introduced in \textbf{\cite{Weil79}} the class $\mathcal{K}^n_s(j)$ consisting of all origin-symmetric convex bodies $K \in \mathcal{K}^n$
for which there exists a \emph{non-negative} Borel measure $\varrho_j(K,\cdot)$ on $\mathrm{Gr}_{j,n}$ such that
\[\mathrm{vol}_j(K|\,\cdot\,) = \mathrm{C}_j \varrho_j(K,\cdot). \]
The classes $\mathcal{K}^n_s(j)$ were subsequently investigated by Goodey and Weil \textbf{\cite{GoodWei91}} and were recently shown to play
an important role in the theory of valuations by Parapatits and Wannerer \textbf{\cite{ParWan13}}. Note that by (\ref{suppvol1}) and (\ref{zonoidchar}), we have $\mathcal{K}^n_s(1) = \mathcal{Z}^n_s$ and that, by (\ref{radsurfcos}), $\mathcal{K}^n_s(n-1)$ coincides with the space $\mathcal{K}^n_s$ of all origin-symmetric convex bodies in $\mathbb{R}^n$. Moreover, a result of Weil \textbf{\cite{Weil79}}
(cf.\ \textbf{\cite[\textnormal{Theorem 5.3.5}]{Schneider14}}) shows that $\mathcal{Z}^n_s \subseteq \mathcal{K}^n_s(j)$ for every $1 \leq j \leq n - 1$.

A convex body $K \in \mathcal{K}^n$ is said to be of class $C^{\infty}_+$ if the boundary of $K$ is a smooth submanifold of $\mathbb{R}^n$ with everywhere
positive curvature. In this case, each measure $S_j(K,\cdot)$, $0 \leq j \leq n - 1$, is absolutely continuous with respect to the spherical Lebesgue measure and has a strictly positive smooth density $s_j(K,\cdot)$ which is (up to a constant) given by the $j$th elementary symmetric function of the principal radii of curvature of $K$.
On the other hand, if $P \in \mathcal{K}^n$ is a polytope, then
\begin{equation} \label{sjpoly}
S_j(P,\omega) = {n - 1 \choose j}^{-1} \sum_{F \in \mathcal{F}^j(P)} \mathcal{H}^{n-1-j}(N(P,F) \cap \omega)\mathrm{vol}_j(F)
\end{equation}
for every Borel set $\omega \subseteq \mathbb{S}^{n-1}$. Here, $\mathcal{F}^{j}(P)$ is the set of all $j$-dimensional faces of $P$, $\mathcal{H}^m$
is $m$-dimensional Hausdorff measure, and $N(P,F)$ denotes the set of all outer unit normal vectors to $P$ at the points of $F$.
Since $N(P,F)$ lies in an $n-1-j$ dimensional great sphere, it follows that $S_j(P,\cdot)$ is concentrated on the union of
finitely many such subspheres. The following converse of this observation was obtained by Goodey and Schneider.

\begin{prop}[\!\!\cite{GoodeySchneider}] \label{sjpolytope} Suppose that $1 \leq j \leq n - 1$ and let $K \in \mathcal{K}^n$ with $\dim K \geq j + 1$.
If the support of the area measure $S_j(K,\cdot)$ can be covered by finitely many $n-1-j$ dimensional great spheres, then $K$ is a polytope.
\end{prop}

The center of mass (centroid) of every area measure of a convex
body is at the origin, that is, for every $1 \leq j \leq n - 1$ and $K \in \mathcal{K}^n$, we have
\begin{equation*}
\int_{\mathbb{S}^{n-1}} u\,dS_j(K,u) = o.
\end{equation*}

The set $\mathcal{S}_j = \{S_j(K,\cdot): K \in \mathcal{K}^n\}$ is dense in the set $\mathcal{M}_o^+(\mathbb{S}^{n-1})$ of all non-negative Borel measures on $\mathbb{S}^{n-1}$ with centroid at
the origin if and only if $j = n-1$. However, $\mathcal{S}_j - \mathcal{S}_j$, $1 \leq j \leq n - 1$, is dense in the
space $\mathcal{M}_o(\mathbb{S}^{n-1})$ of \emph{signed} Borel measures on $\mathbb{S}^{n-1}$ with centroid at the origin
(cf.\ \textbf{\cite[\textnormal{p.\ 477}]{Schneider14}}).

The classical \emph{Minkowski problem} asks for necessary and sufficient conditions for a Borel measure on $\mathbb{S}^{n-1}$ to be
the surface area measure of a convex body. Its solution is one of the cornerstones of the Brunn-Minkowski theory and known as \emph{Minkowski's existence theorem}
(see, e.g., \textbf{\cite[\textnormal{Theorem 8.2.2}]{Schneider14}}): A non-negative Borel measure $\mu$ on $\mathbb{S}^{n-1}$ is the
surface area measure of some convex body in $\mathbb{R}^n$ if and only if $\mu$ is not concentrated on any great subsphere of $\mathbb{S}^{n-1}$ and has its centroid at the origin.

The analogue of Minkowski's problem for the first-order area measure is known as the \emph{Christoffel problem}.
In order to describe its solution by Berg \textbf{\cite{cberg}}, recall that, for $K \in \mathcal{K}^n$, the measure $S_1(K,\cdot)$
and the support function $h(K,\cdot)$ are related by a \emph{linear} differential operator $\Box_n$ in the following way
\begin{equation} \label{boxhks1}
S_1(K,\cdot) = h(K,\cdot) +  \frac{1}{n-1}\Delta_{\mathbb{S}}h(K,\cdot) =: \Box_n h(K,\cdot).
\end{equation}
Here, $\Delta_{\mathbb{S}}$ denotes the Laplacian on $\mathbb{S}^{n-1}$ and equation (\ref{boxhks1}) has to be understood in the sense of distributions.
Since $\Delta_{\mathbb{S}} Y_k = -k(k + n - 2)\,Y_k$ for every $Y_k \in \mathcal{H}_k^n$, the definition of $\Box_n$ implies that, for $f \in C^{\infty}(\mathbb{S}^{n-1})$, the spherical harmonic
expansion of $\Box_n f$ is given by
\begin{equation} \label{boxnmult}
\Box_n f \sim \sum_{k=0}^{\infty} \frac{(1-k)(k+n-1)}{n-1} \pi_kf.
\end{equation}
In particular, the kernel of $\Box_n$ is the space $\mathcal{H}_1^n$ consisting of the restrictions of linear functions on $\mathbb{R}^n$ to $\mathbb{S}^{n-1}$. If we let
\[C_{o}^{\infty}(\mathbb{S}^{n-1}) := \{f \in C^{\infty}(\mathbb{S}^{n-1}): \pi_1 f = 0\},  \]
then $\Box_n: C_{o}^{\infty}(\mathbb{S}^{n-1}) \rightarrow C_{o}^{\infty}(\mathbb{S}^{n-1})$ is an $\mathrm{SO}(n)$ equivariant isomorphism of topological vector spaces.
An explicit expression for the inverse of $\Box_n$ was obtained by Berg \textbf{\cite{cberg}}. He proved that for every $n \geq 2$ there exists a uniquely
determined $C^{\infty}$ function $g_n$ on $(-1,1)$ such that the associated zonal function $\breve{g}_n(u)=g_n(u\cdot \bar{e})$ is in $L^1(\mathbb{S}^{n-1})$ and
\begin{equation} \label{multgn}
a_1^n[g_n] = 0, \quad \qquad a_k^n[g_n] = \frac{n-1}{(1-k)(k+n-1)}, \quad k \neq 1.
\end{equation}
It follows from (\ref{boxnmult}), (\ref{funkheckgen}), and (\ref{multgn}) that for every $f \in C_\mathrm{o}^{\infty}(\mathbb{S}^{n-1})$,
\begin{equation} \label{boxninverse}
f = (\Box_n f) \ast \breve{g}_n.
\end{equation}
From (\ref{boxninverse}), Berg concluded that a measure $\mu \in \mathcal{M}_o^+(\mathbb{S}^{n-1})$ is the first-order area measure of a
convex body in $\mathbb{R}^n$ if and only if $\mu \ast \breve{g}_n$ is a support function.

At the end of this section, we need the following generalization of (\ref{boxninverse}) that follows from a recent result of Goodey and Weil \textbf{\cite[\textnormal{Theorem 4.3}]{GooWei14}}:
For every $j \in \{2, \ldots, n\}$, the convolution transform
\[\mathrm{T}_{g_j}: C_{\mathrm{o}}^{\infty}(\mathbb{S}^{n-1}) \rightarrow C_{\mathrm{o}}^{\infty}(\mathbb{S}^{n-1}), \quad f \mapsto f \ast \breve{g}_j,  \]
is an isomorphism. Let $\Box_j: C_{\mathrm{o}}^{\infty}(\mathbb{S}^{n-1}) \rightarrow C_{\mathrm{o}}^{\infty}(\mathbb{S}^{n-1})$ denote its inverse.

\vspace{0.2cm}

The problem of finding necessary and sufficient conditions for a Borel measure on $\mathbb{S}^{n-1}$ to be an intermediate area measure
of a convex body is known as the \emph{Christoffel-Minkowski problem} and has only been partially resolved (see, e.g., \linebreak \textbf{\cite[\textnormal{Chapter 8.4}]{Schneider14}}).
An in-depth analysis of the problem under additional regularity assumptions was carried out by Guan et al. \textbf{\cite{GuaMa03, GuaLinMa, GuaMaZhou}}.
The following corollary to one of their results \textbf{\cite[\textnormal{Theorem 1.3}]{GuaMa03}} is of particular interest to us.

\begin{prop} \label{guanma} Suppose that $1 \leq j \leq n - 1$. If $K \in \mathcal{K}^n_s$ is an origin-symmetric convex body of class $C^{\infty}_+$,
then $\rho(K,\cdot)^j \in C_e^{\infty}(\mathbb{S}^{n-1})$ is the density of the area measure of order $j$ of a convex body $L \in \mathcal{K}^n_s$ with non-empty interior.
\end{prop}

Firey \textbf{\cite{firey70}} gave the following solution of the Christoffel-Minkowski problem for sufficiently regular \emph{convex bodies of revolution}.
When considering such bodies, we will always assume that they are $\mathrm{SO}(n-1)$ invariant, that is, their axes of revolution is the line spanned by $\bar{e} \in \mathbb{S}^{n-1}$.

\begin{theorem}[\!\!\cite{firey70}] \label{firey_area_measures}Suppose that $1 \leq j \leq n - 2$. A zonal function $s(\bar{e}\cdot\,.\,)$ on $\mathbb{S}^{n-1}$ is the density of
the area measure $S_j(K,\cdot)$ of a (sufficiently smooth) strictly convex body of revolution $K \in \mathcal{K}^n$ if and only if $s$ satisfies the following conditions:
\begin{enumerate}
\item[(i)] $s$ is continuous on $(-1,1)$ and $\lim_{t\rightarrow \pm 1} s(t)$ is finite;
\item[(ii)] $\int_t^1 \xi\,s(\xi)(1-\xi^2)^{\frac{n-3}{2}}d\xi>0$ for $t \in (-1,1)$ and vanishes for $t = -1$;
\item[(iii)] $s(t)(1-t^2)^{\frac{n-1}{2}} > (n-1-j)\int_t^1 \xi\,s(\xi)(1-\xi^2)^{\frac{n-3}{2}}d\xi$ for all $t \in (-1,1)$.
\end{enumerate}
\end{theorem}

\vspace{0.1cm}

We turn now to the theory of valuations. Let $\mathbf{Val}$ denote the vector space of continuous, translation invariant,
scalar-valued valuations. The basic structural result about  $\mathbf{Val}$ is McMullen's decomposition theorem (cf.\ \textbf{\cite[\textnormal{Theorem 6.3.1}]{Schneider14}}):
\begin{equation} \label{mcmullen}
\mathbf{Val} = \bigoplus_{0\leq j \leq n} \left ( \mathbf{Val}_j^+ \oplus \mathbf{Val}_j^-\right ).
\end{equation}
Here, $\mathbf{Val}_j^\pm \subseteq \mathbf{Val}$ denote the subspaces of \emph{even/odd} valuations (homogeneous) of degree $j$.
Using (\ref{mcmullen}), it is not difficult to show that the space $\mathbf{Val}$ becomes a Banach space, when endowed with
the norm $\|\phi \|=\sup \{|\phi(K)|: K \subseteq B\}$. The natural continuous action of the general linear group $\mathrm{GL}(n)$ on this Banach
space is defined as follows: For $A \in \mathrm{GL}(n)$ and every $K \in \mathcal{K}^n$,
\[(A\phi)(K)=\phi(A^{-1}K), \qquad \phi \in \mathbf{Val}.  \]

Suppose that $1 \leq j \leq n - 1$ and let $\mathcal{M}(\mathrm{Gr}_{j,n})$ denote the space of \emph{signed} Borel measures on $\mathrm{Gr}_{j,n}$.
By the Irreducibility Theorem of Alesker \textbf{\cite{Alesker01}}, the map $\mathrm{Cr}_j: \mathcal{M}(\mathrm{Gr}_{j,n}) \rightarrow \mathbf{Val}_j^+$, defined by
\[(\mathrm{Cr}_j\mu)(K) = \int_{\mathrm{Gr}_{j,n}} \mathrm{vol}_j(K|E)\,d\mu(E),   \]
has dense image. This motivates the following notion.

\vspace{0.3cm}

\noindent {\bf Definition} \emph{A measure $\mu \in \mathcal{M}(\mathrm{Gr}_{j,n})$, $1 \leq j \leq n - 1$, is called a
\emph{Crofton measure} for the valuation $\phi \in \mathbf{Val}_j^+$ if $\mathrm{Cr}_j\mu = \phi$.}

\vspace{0.3cm}

In order to state a more precise description of valuations admitting a Crofton measure, we also need to recall the notion of smooth valuations.

\vspace{0.3cm}

\noindent {\bf Definition} \emph{A valuation $\phi \in \mathbf{Val}$ is called smooth if the map $\mathrm{GL}(n) \rightarrow \mathbf{Val}$, defined by $A \mapsto A\phi$, is
infinitely differentiable.}

\pagebreak

The vector space $\mathbf{Val}^{\infty}$ of all smooth translation invariant valuations carries a natural
Fr\'echet space topology (see, e.g., \textbf{\cite{SchuWann16}}) which is stronger then the Banach space topology on $\mathbf{Val}$.
Let $\mathbf{Val}_j^{\pm,\infty}$ denote the subspaces of smooth valuations in $\mathbf{Val}_j^{\pm}$. A basic fact from representation theory implies that the spaces of
smooth valuations $\mathbf{Val}_j^{\pm,\infty}$ are $\mathrm{GL}(n)$ invariant dense subspaces of $\mathbf{Val}_j^{\pm}$.

Suppose that $1 \leq j \leq n - 1$. The {\it Klain map}
\[\mathrm{Kl}_j: \mathbf{Val}_j^+ \rightarrow C(\mathrm{Gr}_{j,n}), \qquad \phi \mapsto \mathrm{Kl}_j\phi,\]
is defined as follows: For $\phi \in \mathbf{Val}_j^+$ and every $E \in \mathrm{Gr}_{j,n}$, consider the restriction $\phi|_E$ of $\phi$
to convex bodies in $E$. This is a continuous, translation invariant valuation of degree $j$ in $E$. Therefore, a classical result of Hadwiger
(see, e.g., \textbf{\cite[\textnormal{Theorem 6.4.8}]{Schneider14}}) implies that $\phi|_E = (\mathrm{Kl}_j\phi)(E)\,\mathrm{vol}_j$,
where $(\mathrm{Kl}_j\phi)(E)$ is a constant depending only on $E$.
The continuous function $\mathrm{Kl}_j\phi \in C(\mathrm{Gr}_{j,n})$ defined in this way is called the \emph{Klain function} of $\phi$.
It is not difficult to see that the map $\mathrm{Kl}_j$ is $\mathrm{SO}(n)$ equivariant and that smooth valuations are mapped to smooth functions, that is,
$\mathrm{Kl}_j: \mathbf{Val}_j^{+,\infty} \rightarrow C^{\infty}(\mathrm{Gr}_{j,n})$.
Moreover, an important result of Klain \textbf{\cite{Klain00}} states that the Klain map $\mathrm{Kl}_j$ is injective for every $j \in \{1, \ldots, n - 1\}$.

Let us now consider the restriction of the Crofton map $\mathrm{Cr}_j$, $1 \leq j \leq n - 1$, to smooth functions.
It is not difficult to see that $\mathrm{Cr}_jf \in \mathbf{Val}_j^{+,\infty}$ for every $f \in C^{\infty}(\mathrm{Gr}_{j,n})$.
Moreover, the Klain function of $\mathrm{Cr}_jf$ is equal to the cosine transform $\mathrm{C}_jf$ of $f$, that is,
\begin{equation} \label{klaincrofton}
\mathrm{Kl}_j \circ \mathrm{Cr}_j = \mathrm{C}_j.
\end{equation}
From this and the main result of \textbf{\cite{AlBern}}, Alesker \textbf{\cite[\textnormal{p.\ 73}]{Alesker03}} deduced the
following:

\begin{theorem}[\!\!\cite{AlBern, Alesker03}] \label{klaincoscrof} Let $1 \leq j \leq n - 1$. The image of the Klain map \linebreak
$\mathrm{Kl}_j: \mathbf{Val}_j^{+,\infty} \rightarrow C^{\infty}(\mathrm{Gr}_{j,n})$ coincides with the image of the cosine transform
$\mathrm{C}_j: C^{\infty}(\mathrm{Gr}_{j,n}) \rightarrow C^{\infty}(\mathrm{Gr}_{j,n})$. Moreover, every smooth valuation $\phi \in \mathbf{Val}_j^{+,\infty}$ admits
a (not necessarily unique for $2 \leq j \leq n - 2$) smooth Crofton measure.
\end{theorem}

Next, we recall the definition of the Alesker-Fourier transform
\[\mathbb{F}: \mathbf{Val}^{+,\infty}_j \rightarrow \mathbf{Val}^{+,\infty}_{n-j}, \qquad 1 \leq j \leq n - 1,  \]
of even valuations (for the odd case, which is much more involved and will not be needed in this article, see \textbf{\cite{Alesker11}}): If
$\phi \in \mathbf{Val}^{+,\infty}_j$, then $\mathbb{F}\phi \in \mathbf{Val}^{+,\infty}_{n-j}$ is the valuation with Klain function given by
\begin{equation} \label{deffourier}
\mathrm{Kl}_{n-j}(\mathbb{F}\phi) = (\mathrm{Kl}_j\phi)^{\bot}.
\end{equation}
By (\ref{cifbot}) and Theorem \ref{klaincoscrof}, the map $\mathbb{F}$ is a well defined $\mathrm{SO}(n)$ equivariant involution.
Moreover, (\ref{klaincrofton}) implies that if $\mu \in \mathcal{M}(\mathrm{Gr}_{j,n})$ is a (smooth) Crofton measure of $\phi$, then $\mu^{\bot} \in \mathcal{M}(\mathrm{Gr}_{n-j,n})$ is
a Crofton measure of $\mathbb{F}\phi$.

In order to define the notion of spherical valuations,
let us first recall the decomposition of the subspaces $\mathbf{Val}_j$ and $\mathbf{Val}^{\infty}_j$ of $j$-homogeneous valuations into $\mathrm{SO}(n)$ irreducible subspaces.

\pagebreak

\begin{theorem}[\!\!\cite{ABS11}] \label{thm_decomposition}
For $0 \leq j \leq n$, the spaces $\mathbf{Val}_j$ and $\mathbf{Val}_j^{\infty}$ are multiplicity free under the action of $\mathrm{SO}(n)$.
Moreover, the highest weights of the $\mathrm{SO}(n)$ irreducible subspaces in either of them are given by the tuples $(\lambda_1,\ldots,\lambda_{\lfloor n/2 \rfloor})$ satisfying (\ref{heiwei}) and the following additional conditions:
\[(i)\ \lambda_k = 0 \mbox{ for } k > \min\{j,n-j\}; \quad \, (ii)\ |\lambda_k| \neq 1 \mbox{ for } 1 \leq k \leq \lfloor \mbox{$\frac{n}{2}$} \rfloor; \quad\, (iii)\  |\lambda_2| \leq 2.  \]
\end{theorem}

\vspace{0.1cm}

The notion of spherical representations with respect to $\mathrm{SO}(n-1)$ (compare Section 2) motivates the following.

\vspace{0.3cm}

\noindent {\bf Definition} \emph{For $0 \leq j \leq n$, the spaces $\mathbf{Val}_j^{\mathrm{sph}}$ and $\mathbf{Val}_j^{\infty,\mathrm{sph}}$
of translation invariant, continuous and smooth \emph{spherical valuations} of degree $j$ are defined as the closures (w.r.t.\ the respective topologies) of the direct sum of all
$\mathrm{SO}(n)$ irreducible subspaces in $\mathbf{Val}_j$ and $\mathbf{Val}_j^{\infty}$, respectively,
which are spherical with respect to $\mathrm{SO}(n - 1)$.}

\vspace{0.3cm}

Theorems \ref{thm_decomposition} and \ref{thmspher} imply that $\mathbf{Val}_j^{\mathrm{sph}}$ and $\mathbf{Val}_j^{\infty,\mathrm{sph}}$
are the closures of the direct sum of all $\mathrm{SO}(n)$ irreducible subspaces in $\mathbf{Val}_j$ and $\mathbf{Val}_j^{\infty}$, respectively,
with highest weights $(k,0,\ldots,0)$, $k \in \mathbb{N}$. In particular, by Theorem \ref{thm_decomposition}, we have
\[\mathbf{Val}_1^{(\infty)} = \mathbf{Val}_1^{(\infty),\mathrm{sph}} \qquad \mbox{and} \qquad \mathbf{Val}_{n-1}^{(\infty)} = \mathbf{Val}_{n-1}^{(\infty),\mathrm{sph}}\]
and, by Theorem \ref{thmspher} (b), every $\mathrm{SO}(n - 1)$-invariant valuation in $\mathbf{Val}_j$ or $\mathbf{Val}_j^{\infty}$, $0 \leq j \leq n$, is spherical.
Moreover, the following alternative description of smooth spherical valuations was established in \textbf{\cite{SchuWann16}}.

\begin{prop} \label{smoothsphiso} For $1 \leq j \leq n - 1$, the map
\begin{equation} \label{defei}
\mathrm{E}_j: C_o^{\infty}(\mathbb{S}^{n-1}) \rightarrow \mathbf{Val}_j^{\infty,\mathrm{sph}}, \quad (\mathrm{E}_jf)(K) = \int_{\mathbb{S}^{n-1}} f(u)\,dS_j(K,u),
\end{equation}
is an $\mathrm{SO}(n)$ equivariant isomorphism of topological vector spaces.
\end{prop}

Proposition \ref{smoothsphiso} and a recent result of Bernig and Hug \textbf{\cite[\textnormal{Lemma 4.8}]{BerHug15+}} now imply the following
relation between the Alesker-Fourier transform of \emph{even} spherical valuations and certain Radon transforms of spherical functions.

\begin{prop} \label{aleskfourradonprop} Suppose that $1 \leq j \leq n - 1$. If the \emph{even} spherical valuation $\phi \in \mathbf{Val}_j^{\infty,\mathrm{sph}}$ is given by
$\phi = \mathrm{E}_jf$ for $f \in C_e^{\infty}(\mathbb{S}^{n-1})$, then $\mathbb{F}\phi \in \mathbf{Val}_{n-j}^{\infty,\mathrm{sph}}$ and,
for every $K \in \mathcal{K}^n$,
\begin{equation} \label{aleskfourradon}
(\mathbb{F}\phi)(K) = \frac{\kappa_{n-j}}{\kappa_{j}}\int_{\mathbb{S}^{n-1}}(\mathrm{R}_{1,j}^{-1} \circ \bot_* \circ \mathrm{R}_{1,n-j})f(u)\,dS_{n-j}(K,u).
\end{equation}
\end{prop}

Note that the function under the integral in (\ref{aleskfourradon}) is well defined, since $\mathrm{R}_{1,n-j}f \in C^{\infty}(\mathrm{Gr}_{n-j,n})^{\mathrm{sph}}$ for every $f\! \in\! C_e^{\infty}(\mathbb{S}^{n-1})$ and the Radon transform $\mathrm{R}_{1,j}: C^{\infty}(\mathrm{Gr}_{1,n}) \rightarrow C^{\infty}(\mathrm{Gr}_{j,n})^{\mathrm{sph}}$ is bijective (cf.\ Example \ref{exptransforms} (b)).

Comparing Propositions \ref{spherfourradon} and \ref{aleskfourradonprop}, we obtain the following critical relation between the
Alesker-Fourier transform of spherical valuations and the spherical Fourier transform.

\pagebreak

\begin{coro} \label{aleskfoursphfour} Suppose that $1 \leq j \leq n - 1$. If the \emph{even} spherical valuation $\phi \in \mathbf{Val}_j^{\infty,\mathrm{sph}}$ is given by
$\phi = \mathrm{E}_jf$ for $f \in C_e^{\infty}(\mathbb{S}^{n-1})$, then, for every $K \in \mathcal{K}^n$,
\[(\mathbb{F}\phi)(K) = \frac{j}{(2\pi)^j(n-j)}\int_{\mathbb{S}^{n-1}}(\mathbf{F}_{j-n}f)(u)\,dS_{n-j}(K,u).\]
\end{coro}

\vspace{0.1cm}

From the computation of the multipliers of the Alesker-Fourier transform of spherical valuations in \textbf{\cite{BerHug15+}} and the spherical Fourier transform
in \textbf{\cite{GooYasYas11}}, it follows that Corollary \ref{aleskfoursphfour} also holds without the assumption on the parity.

\vspace{0.1cm}

In the final part of this section, we recall results on Minkowski valuations intertwining rigid motions.
To this end, let $\mathbf{MVal}_j$, $0 \leq j \leq n$, denote the set of all continuous and translation invariant
Minkowski valuations of degree~$j$. It is easy to see (cf.\ \textbf{\cite{Schuster10}}) that if $\Phi_j \in \mathbf{MVal}_j$ is $\mathrm{SO}(n)$ equivariant, then
the $\mathrm{SO}(n-1)$ invariant real valued valuation $\psi_{\Phi_j} \in \mathbf{Val}_j$, defined by
\[\psi_{\Phi_j}(K) = h(\Phi_jK,\bar{e}), \]
uniquely determines $\Phi_j$ and is called \emph{the associated real valued valuation of} $\Phi_j$.
Motivated by this simple fact, the following definition was first given in \textbf{\cite{Schuster10}}.

\vspace{0.3cm}

\noindent {\bf Definition} \emph{An $\mathrm{SO}(n)$ equivariant Minkowski valuation $\Phi_j \in \mathbf{MVal}_j$, $0 \leq j \leq n$, is called \emph{smooth}
if its associated real valued valuation $\psi_{\Phi_j} \in \mathbf{Val}_j$ is smooth.}

\vspace{0.3cm}

It was recently proved in \textbf{\cite{SchuWann16}} that
any $\mathrm{SO}(n)$ equivariant Minkowski valuation in $\mathbf{MVal}_j$, $0 \leq j \leq n$, can be approximated by smooth ones.

The systematic investigation of Minkowski valuations intertwining affine transformations has its origin in two articles by Ludwig \textbf{\cite{Ludwig02, Ludwig05}} and has
become a research focus in integral geometry in recent years (see, e.g., \textbf{\cite{AbardiaBernig11, Hab12, SchuWann12, Wannerer11}}).
The classification problem of Minkowski valuations intertwining rigid motions only has proven to be more difficult and
has not been completely resolved. The following theorem contains two results of Wannerer and the second author on convolution representations of such Minkowski valuations
which are even.

\begin{theorem}[\!\!\cite{SchuWann15, SchuWann16}] \label{repschuwan} Suppose that $1 \leq j \leq n - 1$ and let
$\Phi_j \in \mathbf{MVal}_j$ be $\mathrm{SO}(n)$ equivariant and even.
\begin{enumerate}
\item[(i)] There exists a uniquely determined zonal measure $\mu_{\Phi_j} \in \mathcal{M}_e(\mathbb{S}^{n-1})$, called the \emph{generating measure} of $\Phi_j$,
such that for every $K \in \mathcal{K}^n$,
\begin{equation}
h(\Phi_jK,\cdot) = S_j(K,\cdot) \ast \mu_{\Phi_j}.
\end{equation}
If $\Phi_j$ is smooth, then $\mu_{\Phi_j}$ has a smooth density.
\item[(ii)] If $\Phi_j$ is smooth, then there exists a uniquely determined $\mathrm{O}(j) \times \mathrm{O}(n-j)$ invariant smooth measure $\sigma_{\Phi_j} \in \mathcal{M}_e(\mathbb{S}^{n-1})$,
called the \emph{spherical Crofton measure} of $\Phi_j$, such that for every $K \in \mathcal{K}^n$,
\begin{equation}
h(\Phi_jK,\cdot) = \mathrm{vol}_j(K|\,\cdot\,) \ast \sigma_{\Phi_j}.
\end{equation}
\end{enumerate}
\end{theorem}

\pagebreak

We now use associated real valued valuations to extend
the Alesker-Fourier transform (at least partially) to Minkowski valuations.

\vspace{0.3cm}

\noindent {\bf Definition} \emph{Let $\Phi_j \in \mathbf{MVal}_j$ and $\Psi_{n-j} \in \mathbf{MVal}_{n-j}$, $1 \leq j \leq n - 1$, both be $\mathrm{SO}(n)$ equivariant and even. We write $\Psi_{n-j} = \mathbb{F}\Phi_j$ and say $\Psi_{n-j}$ is the Alesker-Fourier transform of $\Phi_j$ if}
\[\mathrm{Kl}_{n-j}(\psi_{\Psi_{n-j}}) = (\mathrm{Kl}_j\psi_{\Phi_j})^{\bot}.\]

\vspace{0.3cm}

Note that if $\Phi_j$ and $\Psi_{n-j}$ are in addition smooth, then $\psi_{\Psi_{n-j}} = \mathbb{F}\psi_{\Phi_j}$ by (\ref{deffourier}).
Moreover, in this case, by Theorem \ref{repschuwan}, both $\Phi_j$ and $\Psi_{n-j}$ admit zonal generating functions $g_{\Phi_j}, g_{\Psi_{n-j}} \in C^{\infty}_e(\mathbb{S}^{n-1})$ \emph{and} smooth spherical Crofton measures $\sigma_{\Phi_j}, \sigma_{\Psi_{n-j}} \in \mathcal{M}_e(\mathbb{S}^{n-1})$.
Hence, from the definition of the convolution, it is not difficult to show (cf.\ \textbf{\cite{SchuWann15, SchuWann16}})  that
\[\psi_{\Phi_j} = \mathrm{E}_jg_{\Phi_j} = \mathrm{Cr}_j \widehat{\sigma_{\Phi_j}} \qquad \mbox{and} \qquad \psi_{\Psi_{n-j}} = \mathrm{E}_{n-j}g_{\Psi_{n-j}} = \mathrm{Cr}_{n-j} \widehat{\sigma_{\Psi_{n-j}}},\]
where
\[\widehat{\sigma_{\Phi_j}} := \pi_{j*}\widehat{\pi^*\sigma_{\Phi_j}} \in \mathcal{M}(\mathrm{Gr}_{j,n})  \]
and $\widehat{\sigma_{\Psi_{n-j}}} \in \mathcal{M}(\mathrm{Gr}_{n-j,n})$ is defined similarly. Here, $\pi: \mathrm{SO}(n) \rightarrow \mathbb{S}^{n-1}$ and $\pi_j: \mathrm{SO}(n) \rightarrow \mathrm{Gr}_{j,n}$ denote the canonical projections. Note that this is well defined by the invariance of $\sigma_{\Phi_j}$ and $\sigma_{\Psi_{n-j}}$.
Therefore, letting
\[\sigma_{\Phi_j}^{\bot} :=  \pi_*\widehat{\pi_{n-j}^*\widehat{\sigma_{\Phi_j}}^{\bot}} \in \mathcal{M}_e(\mathbb{S}^{n-1}), \]
Corollary \ref{aleskfoursphfour}, the remark after (\ref{deffourier}), and a standard approximation argument imply the following.

\begin{prop} \label{minkvalfourier} Let $\Phi_j \in \mathbf{MVal}_j$, $1 \leq j \leq n - 1$, be $\mathrm{SO}(n)$ equivariant and even with generating measure $\mu_{\Phi_j} \in \mathcal{M}_e(\mathbb{S}^{n-1})$
and suppose that $\Phi_j$ admits a spherical Crofton measure $\sigma_{\Phi_j} \in \mathcal{M}_e(\mathbb{S}^{n-1})$.
Then $\Psi_{n-j} \in \mathbf{MVal}_{n-j}$ is the Alesker-Fourier transform of $\Phi_j$ if and only if $\frac{j}{(2\pi)^j(n-j)}\mathbf{F}_{j-n}\mu_{\Phi_j}$ and $\sigma_{\Phi_j}^{\bot}$
are the generating measure and spherical Crofton measure of $\Psi_{n-j}$, respectively.
\end{prop}

Note that it is an open problem whether the Alesker-Fourier transform of every smooth even Minkowski valuation which is $\mathrm{SO}(n)$ equivariant and translation invariant
is well defined (cf.\ \textbf{\cite{SchuWann15}}). However, in the following example we exhibit for every $1 \leq j \leq n - 1$ a pair of Minkowski valuations which are related
via the Alesker-Fourier transform.

\begin{expls} \label{exp17} \end{expls}

\vspace{-0.2cm}

\begin{enumerate}
\item[(a)] For $1 \leq j \leq n - 1$, let $\Pi_j \in \mathbf{MVal}_j$ denote the projection body map of order $j$, defined by
\[h(\Pi_jK,u) = V_j(K|u^{\bot}) = \frac{1}{2} \int_{\mathbb{S}^{n-1}} |u\cdot v|\,dS_j(K,v), \qquad u \in \mathbb{S}^{n-1}.   \]

\pagebreak

Note that each $\Pi_j$ is $\mathrm{SO}(n)$ equivariant and even but {\it not} smooth. Moreover, each $\Pi_j$ is injective on origin-symmetric convex bodies with non-empty interior. Their
continuous generating function is given by $g_{\Pi_j}(u) = \frac{1}{2}|\bar{e}\cdot u|$. It follows from (\ref{radsurfcos}) and (\ref{radconv}) that
\begin{equation} \label{pijradon}
h(\Pi_j K,\cdot) = \frac{\kappa_{n-1}}{\kappa_j}\mathrm{R}_{n-j,1}\mathrm{vol}_j(K|\,\cdot\,)^{\bot} = \frac{\kappa_{n-1}}{\kappa_j}\,\mathrm{vol}_j(K|\,\cdot\,) \ast \widehat{\lambda}_{1,n-j}^{\bot}.  \end{equation}
Thus, the measure $\frac{\kappa_{n-1}}{\kappa_j}\widehat{\lambda}_{1,n-j}^{\bot}$ is the spherical Crofton measure of $\Pi_j$.

\item[(b)] For $1 \leq j \leq n - 1$, Goodey and Weil introduced and investigated in \textbf{\cite{GooWei92,GooWei12,GooWei14}} the
{\it normalized} mean section operator of order $j$, denoted by \linebreak $\mathrm{M}_j \in \mathbf{MVal}_{n+1-j}$.
In \textbf{\cite[\textnormal{Theorem 4.4}]{GooWei14}} they proved that
\begin{equation} \label{genfctmi}
h(\mathrm{M}_jK,\cdot) = q_{n,j}\,S_{n+1-j}(K,\cdot) \ast \breve{g}_j,
\end{equation}
where
\begin{equation*} \label{qnj}
q_{n,j} = \frac{j-1}{2\pi(n+1-j)}\,\frac{\kappa_{j-1}\kappa_{j-2}\kappa_{n-j}}{\kappa_{j-3}\kappa_{n-2}}.
\end{equation*}
Hence, the multiple $q_{n,j}\breve{g}_j$ of the Berg function is the generating function of $\mathrm{M}_j$.
Note that $\mathrm{M}_j$ is continuous and $\mathrm{SO}(n)$ equivariant but \emph{not} even. Moreover, $\mathrm{M}_j$ determines a convex body with non-empty interior up to translations.
For the even part of $\mathrm{M}_j$, Goodey and Weil \textbf{\cite{GooWei92}} proved that
\begin{equation} \label{Mjradon}
h(\mathrm{M}_j^+K,\cdot)=\frac{j\kappa_j\kappa_{n-1}}{2n\kappa_{j-1}\kappa_n}\,\mathrm{R}_{n+1-j,1}\mathrm{vol}_{n+1-j}(K|\,\cdot\,).
\end{equation}
Thus, by (\ref{radconv}), $\frac{j\kappa_j\kappa_{n-1}}{2n\kappa_{j-1}\kappa_n}\widehat{\lambda}_{1,n+1-j}$ is the spherical Crofton measure of $\mathrm{M}_j^+$.
Comparing this with Example \ref{exp17} (a), it follows from Proposition \ref{minkvalfourier} that the \emph{renormalized even} mean section operator
\[\overline{\mathrm{M}}_{n-j} := \frac{2n\kappa_n}{(j+1)\kappa_{j+1}}\mathrm{M}_{j+1}^+ \in \mathbf{MVal}_{n-j} \]
is the Alesker-Fourier transform of $\Pi_j$, that is,
\begin{equation} \label{meanfourproj}
\overline{\mathrm{M}}_{n-j} = \mathbb{F}\Pi_j.
\end{equation}
\end{enumerate}

Comparing generating functions of $\Pi_j$ and $\overline{\mathrm{M}}_{n-j}$ and using Proposition \ref{minkvalfourier} as well as (\ref{sphfourinv}),
yields the following representation of the spherical Fourier transform.

\begin{coro} \label{sphfourconvrep} Suppose that $1 \leq j \leq n - 1$. Then
\[\mathbf{F}_{-j} = \frac{(2\pi)^{n-j}j(j+1)\kappa_{j+1}}{4(n-j)n\kappa_n}q_{n,j+1}^{-1} \mathrm{C}_1 \circ \Box_{j+1}.\]
\end{coro}

\pagebreak

\centerline{\large{\bf{ \setcounter{abschnitt}{4}
\arabic{abschnitt}. Proof of the Main Results}}}

\reseteqn \alpheqn \setcounter{theorem}{0}

\vspace{0.6cm}

After these preparations we are now in a position to present the proofs of Theorems \ref{mainthm1} and \ref{mainthm2} as well as Corollary \ref{korotothm3} in this section.

\vspace{0.2cm}

\noindent {\it Proof of Theorem \ref{mainthm1}}. First note that, by Proposition \ref{goodeyaleks} and (\ref{rijbot}), $K$ is the $j$-projection body of $L$ if and only if
\begin{equation} \label{prfthm2a}
(\mathrm{R}_{j,n-1}\mathrm{vol}_j(K|\,\cdot\,))^{\perp} = \mathrm{R}_{n-j,1} \mathrm{vol}_{n-j}(L|\,\cdot\,).
\end{equation}
Using (\ref{radsurfcos}) as well as the fact that $\mathrm{R}_{1,n-j}$ is the adjoint of $\mathrm{R}_{n-j,1}$, we see that (\ref{prfthm2a}) holds if and only if
\begin{equation} \label{prfthm2b}
\frac{\kappa_j}{\kappa_{n-1}}\int_{\mathbb{S}^{n-1}}f(u)V_j(K|u^{\perp})\,du = \int_{\mathrm{Gr}_{n-j,n}}(\mathrm{R}_{1,n-j}f)(E)\mathrm{vol}_{n-j}(L|E)\,dE
\end{equation}
for every $f \in C^{\infty}_e(\mathbb{S}^{n-1})$. Since $\mathrm{R}_{1,n-j}f \in C^{\infty}(\mathrm{Gr}_{n-j,n})^{\mathrm{sph}}$ and
\[\mathrm{R}_{n-j,n-1}: C^{\infty}(\mathrm{Gr}_{n-j,n})^{\mathrm{sph}} \rightarrow C^{\infty}(\mathrm{Gr}_{n-1,n})^{\mathrm{sph}} \]
is bijective, it follows that the integral on the right hand side of (\ref{prfthm2b}) is equal to
\[\int_{\mathrm{Gr}_{n-1,n}}(\mathrm{R}_{n-1,n-j}^{-1}\mathrm{R}_{1,n-j}f)(F)(\mathrm{R}_{n-j,n-1}\mathrm{vol}_{n-j}(L|\,\cdot\,))(F)\,dF.    \]
Since $\bot_*$ is clearly self-adjoint, this can be further rewritten, by using again (\ref{radsurfcos}) and (\ref{rijbot}), to obtain
\[\frac{\kappa_{n-j}}{\kappa_{n-1}}\int_{\mathbb{S}^{n-1}}(\mathrm{R}_{1,j}^{-1}\circ \bot_* \circ \mathrm{R}_{1,n-j})f(u)V_{n-j}(L|u^{\perp}\,)\,du.    \]
Consequently, by Proposition \ref{spherfourradon}, (\ref{prfthm2a}) holds if and only if
\[\int_{\mathbb{S}^{n-1}}f(u)V_j(K|u^{\perp})\,du = \frac{j}{(2\pi)^j(n-j)}\int_{\mathbb{S}^{n-1}}(\mathbf{F}_{j-n}f)(u)V_{n-j}(L|u^{\perp}\,)\,du\]
for every $f \in C^{\infty}_e(\mathbb{S}^{n-1})$. Since the cosine transform $\mathrm{C}_1$ is self-adjoint, it follows from (\ref{radsurfcos}) and the obvious fact
that the multiplier transformations $\mathbb{F}_{j-n}$ and $\mathrm{C}_1$ commute that this is equivalent to
\[\int_{\mathbb{S}^{n-1}}\mathrm{C}_1f(u)\,dS_j(K,u) = \frac{j}{(2\pi)^j(n-j)}\int_{\mathbb{S}^{n-1}}(\mathbf{F}_{j-n}\mathrm{C}_1f)(u)\,dS_{n-j}(L,u)\]
for every $f \in C^{\infty}_e(\mathbb{S}^{n-1})$. Substituting $h = \mathbf{F}_{j-n}\mathrm{C}_1f$ and using (\ref{sphfourinv}), we finally obtain the desired relation
\[\int_{\mathbb{S}^{n-1}}\mathbf{F}_{-j}h(u)\,dS_j(K,u)=\frac{(2\pi)^{n-j}j}{(n-j)}\int_{\mathbb{S}^{n-1}}h(u)\,dS_{n-j}(L,u) \]
for every $h \in C^{\infty}_e(\mathbb{S}^{n-1})$ which completes the proof. \hfill $\blacksquare$

\pagebreak

Before we continue, we include here also a short proof of Theorem \ref{koldmain} which underlines the dual nature of Theorems \ref{koldmain} and \ref{mainthm1} and
shows that for the 'if' part of the statement no additional regularity assumptions are required.

\vspace{0.2cm}

\noindent {\it Proof of Theorem \ref{koldmain}}. By passing to polar coordinates, it follows that $D$ is the $j$-intersection body of $M$ if and only if
\begin{equation} \label{prfthm1a}
\frac{\kappa_j}{\kappa_{n-j}}\!\int_{\mathrm{Gr}_{n-j,n}}\!\!\!\!\!\!\!\!f(E)(\mathrm{R}_{1,j}\rho(D,\cdot)^j)^{\perp}(E)\,dE =
\int_{\mathrm{Gr}_{n-j,n}}\!\!\!\!\!\!\!\!f(E) \mathrm{R}_{1,n-j} \rho(M,\cdot)^{n-j}(E)\,dE
\end{equation}
for every $f \in C^{\infty}(\mathrm{Gr}_{n-j,n})$. Since $\bot_*$ is self-adjoint and $\mathrm{R}_{j,i}$ is the adjoint of $\mathrm{R}_{i,j}$,
(\ref{prfthm1a}) is equivalent to
\begin{equation*}
\frac{\kappa_j}{\kappa_{n-j}}\int_{\mathbb{S}^{n-1}}(\mathrm{R}_{j,1}f^{\bot})(u)\rho(D,u)^j\,du = \int_{\mathbb{S}^{n-1}}\mathrm{R}_{n-j,1}f(u)\rho(M,u)^{n-j}\,du
\end{equation*}
for every $f \in C^{\infty}(\mathrm{Gr}_{n-j,n})$. Substituting now $h = \mathrm{R}_{n-j,1}f$ and using that the Radon transform
$\mathrm{R}_{n-j,1}: C^{\infty}(\mathrm{Gr}_{n-j,n}) \rightarrow C^{\infty}_e(\mathbb{S}^{n-1})$ is surjective, we see that $D$ is the $j$-intersection body of $M$ if and only if
\begin{equation} \label{prfthm1b}
\frac{\kappa_j}{\kappa_{n-j}}\int_{\mathbb{S}^{n-1}}(\mathrm{R}_{j,1}\circ \bot_* \circ \mathrm{R}_{n-j,1}^{-1})h(u)\rho(D,u)^j\,du = \int_{\mathbb{S}^{n-1}}h(u)\rho(M,u)^{n-j}\,du
\end{equation}
for all $h \in C^{\infty}_e(\mathbb{S}^{n-1})$. Since $\mathbf{F}_{-j}$ and $\bot_*$ are self-adjoint and $\mathrm{R}_{j,i}$ is the adjoint of $\mathrm{R}_{i,j}$,
it follows from Proposition \ref{spherfourradon} that
\[\mathbf{F}_{-j} = \frac{(2\pi)^{n-j}\,j\,\kappa_j}{(n-j)\,\kappa_{n-j}}\,\mathrm{R}_{j,1}\circ \bot_* \circ \mathrm{R}_{n-j,1}^{-1}.  \]
Hence, (\ref{prfthm1b}) is equivalent to
\[\int_{\mathbb{S}^{n-1}}\mathbf{F}_{-j}h(u)\rho(D,u)^j\,du=\frac{(2\pi)^{n-j}j}{n-j}\int_{\mathbb{S}^{n-1}}h(u)\rho(M,u)^{n-j}\,du \]
for every $h \in C^{\infty}_e(\mathbb{S}^{n-1})$. \hfill $\blacksquare$

\vspace{0.3cm}

With our next result, we complete the proofs of Theorem \ref{mainthm2} and Corollary \ref{korotothm3}.

\begin{theorem} \label{mainequiv} Let $1 \leq j \leq n - 1$ and let $K$ and $L$ be origin-symmetric convex bodies with non-empty interior in $\mathbb{R}^n$.
Then the following statements are equivalent:
\begin{enumerate}
\item[(i)] $K$ is the $j$-projection body of $L$;
\item[(ii)] $\mathbf{F}_{-j}S_j(K,\cdot) = \frac{(2\pi)^{n-j}j}{n-j}S_{n-j}(L,\cdot)$;
\item[(iii)] $\Pi_jK = \overline{\mathrm{M}}_{n-j}L$;
\item[(iv)] $\phi(K) = (\mathbb{F}\phi)(L)$ for all even $\phi \in \mathbf{Val}_j^{\infty,\mathrm{sph}}$.
\end{enumerate}
\end{theorem}

\pagebreak

\noindent {\it Proof}. We have already seen that (i) and (ii) are equivalent. In order to prove that (ii) and (iii) are equivalent, we first note that, by the
definition of $\Pi_j$ and $\overline{\mathrm{M}}_{n-j}$ and (\ref{genfctmi}), (iii) is equivalent to
\[h(\Pi_jK,\cdot)=\frac{1}{2} \mathrm{C}_1 S_j(K,\cdot) = \frac{2n\kappa_n}{(j+1)\kappa_{j+1}} q_{n,j+1}\,S_{n-j}(L,\cdot) \ast \breve{g}_{j+1}=h(\overline{\mathrm{M}}_{n-j}L,\cdot).  \]
Since convolution transforms are self-adjoint, integrating both sides yields that this is equivalent to
\[ \frac{(j+1)\kappa_{j+1}}{4n\kappa_n} q_{n,j+1}^{-1} \int_{\mathbb{S}^{n-1}} \mathrm{C}_1f(u)\,dS_j(K,u) =  \int_{\mathbb{S}^{n-1}} (f \ast \breve{g}_{j+1})(u)\,dS_{n-j}(L,u)  \]
for every $f \in C_e^\infty(\mathbb{S}^{n-1})$. Substituting $h = f \ast \breve{g}_{j+1}$ and using Corollary \ref{sphfourconvrep}, it follows that this holds if and only if
\begin{equation} \label{ii}
\int_{\mathbb{S}^{n-1}}\mathbf{F}_{-j}h(u)\,dS_j(K,u)=\frac{(2\pi)^{n-j}j}{(n-j)}\int_{\mathbb{S}^{n-1}}h(u)\,dS_{n-j}(L,u)
\end{equation}
for every $h \in C^{\infty}_e(\mathbb{S}^{n-1})$ which is precisely (ii).

Finally, in order to see that (iv) is equivalent to (ii), first note that, by Proposition \ref{smoothsphiso} and Corollary \ref{aleskfoursphfour}, (iv) is equivalent to
\[\int_{\mathbb{S}^{n-1}}f(u)\,dS_j(K,u) = \frac{j}{(2\pi)^j(n-j)}\int_{\mathbb{S}^{n-1}}(\mathbf{F}_{j-n}f)(u)\,dS_{n-j}(L,u)\]
for every $f \in C_e^\infty(\mathbb{S}^{n-1})$. Substituting this time $h = \mathbf{F}_{j-n}f$ and using (\ref{sphfourinv}), it follows that this holds if and only if (\ref{ii}) holds
for every $h \in C^{\infty}_e(\mathbb{S}^{n-1})$, which completes the proof.
\hfill $\blacksquare$

\vspace{0.3cm}

We remark that $\Pi_j$ and $\overline{\mathrm{M}}_{n-j}$ can be replaced in statement (iii) by any pair of Minkowski valuations intertwining rigid motions which are
related by the Alesker-Fourier transform and injective on origin-symmetric convex bodies. This follows easily from Proposition \ref{minkvalfourier}.
Next, note that, by Theorem \ref{klaincoscrof}, (\ref{klaincrofton}), and (\ref{deffourier}), the space $\mathbf{Val}_j^{\infty,\mathrm{sph}}$ of smooth spherical valuations can be replaced in (iv) by
the entire space $\mathbf{Val}_j^{\infty}$, leading however to a weaker statement.

Finally, we also note that Corollary \ref{korotothm3} follows also directly from (\ref{pijradon}), (\ref{Mjradon}), and Proposition \ref{goodeyaleks}.
Together with the arguments of the first part of the proof of Theorem \ref{mainequiv}, this can be used to give an alternative proof of Theorem \ref{mainthm1}.

\vspace{1cm}

\centerline{\large{\bf{ \setcounter{abschnitt}{5}
\arabic{abschnitt}. Properties and Examples of $j$-Projection Bodies}}}

\reseteqn \alpheqn \setcounter{theorem}{0}

\vspace{0.6cm}

In this section, we first prove that the class of $j$-projection bodies is invariant under non-degenerate linear transformations.
We then collect several examples of $j$-projection bodies from the literature and compute a family of new examples using Theorem \ref{mainthm1}.
At the end of the section, we generalize two more well known properties of the classical $1$-projection bodies to all $j > 1$.

\pagebreak

The fact that the class of $j$-intersection bodies is invariant under the general linear group $\mathrm{GL}(n)$ was first observed by Milman \textbf{\cite{Milman06}}.
The proof of the following dual counterpart is based on ideas of Schneider \textbf{\cite{Schneider96}}, who proved it in the special case described in Example \ref{exp21} (d) below.

\begin{theorem} \label{thmglninv} Let $1 \leq j \leq n - 1$, $A \in \mathrm{GL}(n)$, and let $K$ and $L$ be origin-symmetric convex bodies with non-empty interior in $\mathbb{R}^n$.
If $K$ is the $j$-projection body of $L$, then $AK$ is the $j$-projection body of
\[|\!\det A|^{\frac{1}{n-j}}\,A^{-\mathrm{T}}L.  \]
In particular, the class of $j$-projection bodies is $\mathrm{GL}(n)$ invariant.
\end{theorem}

\noindent {\it Proof}. First let $|\!\det A| = 1$. Using the polar decomposition of $A$ and the fact that the statement is clearly true
for orthogonal linear maps, we may assume that $A$ is \emph{symmetric} and \emph{positive definite}. Thus, we have to show that
\begin{equation} \label{glninv}
\mathrm{vol}_{j}(AK|E^{\bot}) = \mathrm{vol}_{n-j}(A^{-1}L|E)
\end{equation}
for every $E \in \mathrm{Gr}_{n-j,n}$. To this end, let $F \in \mathrm{Gr}_{j,n}$ and $u_1, \ldots, u_j$ be a set of orthonormal vectors in $F$.
In order to compute $\mathrm{vol}_j(AK|F)$, we may assume, using the singular value decomposition of $A$, that the vectors $Au_1, \ldots, Au_j$ are also orthogonal.
Then,
\begin{eqnarray*}
\mathrm{vol}_j(AK|F) & = & \mathrm{vol}_j\left ( \left \{ \sum_{i=1}^j (Ax \cdot u_i)u_i: x \in K \right \} \right )  \\
 & = & \left ( \prod_{i=1}^j \|Au_i\|_2 \right )\mathrm{vol}_j(K|AF) = c_j(A,F)\,\mathrm{vol}_j(K|AF),
\end{eqnarray*}
where $c_j(A,F)$ depends on $A$ and $F$ only and not on $K$. Consequently, using that $AE^{\bot} = (A^{-1}E)^{\bot}$ for every $E \in \mathrm{Gr}_{n-j,n}$, we obtain
\begin{eqnarray*}
\mathrm{vol}_j(AK|E^{\bot}) & = & c_j(A,E^{\bot})\mathrm{vol}_j(K|(A^{-1}E)^{\bot})= c_j(A,E^{\bot})\mathrm{vol}_{n-j}(L|A^{-1}E) \\
& = & \frac{c_j(A,E^{\bot})}{c_{n-j}(A^{-1},E)}\mathrm{vol}_{n-j}(A^{-1}L|E).
\end{eqnarray*}
Choosing now $K$ to be the unit cube in $\mathbb{R}^n$, it follows from a result of Schnell \textbf{\cite{Schnell94}} (see Example \ref{exp21} (b) below) that
$c_j(A,E^{\bot}) = c_{n-j}(A^{-1},E)$ which yields the desired equation (\ref{glninv}).

Finally, let $A \in \mathrm{GL}(n)$ be arbitrary. Then, by the first part of the proof,
\begin{eqnarray*}
\mathrm{vol}_j(AK|E^{\bot}) \!\!\! & = \!\!\! & |\!\det A|^{\frac{j}{n}}\mathrm{vol}_{j}(|\!\det A|^{-\frac{1}{n}}AK|E^{\bot}) \\
& = \!\!\! & |\!\det A|^{\frac{j}{n}}\mathrm{vol}_{n-j}(|\!\det A|^{\frac{1}{n}}A^{-\mathrm{T}}L|E) = \mathrm{vol}_{n-j}(|\!\det A|^{\frac{1}{n-j}}A^{-\mathrm{T}}L|E)
\end{eqnarray*}
holds for every $E \in \mathrm{Gr}_{n-j,n}$ as desired. \hfill $\blacksquare$

\pagebreak

Theorem \ref{thmglninv} suggests that it should be possible to define the notion of \linebreak $j$-projection bodies
in $\mathrm{SL}(n)$ invariant terms without referring to any Euclidean structure. Indeed, the authors are obliged to
S.\ Alesker for communicating such an $\mathrm{SL}(n)$ invariant definition to us which we will state in the following. It requires
a basic familiarity with the notion of a line bundle over a manifold.

Let $V$ be an $n$-dimensional real vector space with a fixed volume form and let $V^*$ denote its dual space.
For $1 \leq j \leq n - 1$, we write $\mathrm{Gr}_j(V)$ for the Grassmannian of all $j$-dimensional subspaces of $V$ and $\mathcal{K}(V)$ for the space of
convex bodies in $V$. Moreover, if $W$ is a finite dimensional vector space, we denote by $\mathrm{Dens}(W)$ the $1$-dimensional
space of Lebesgue measures on $W$. Finally, for any measurable subset $M$ of $W$, we define an element $\mathrm{ev}_M \in \mathrm{Dens}(W)^*$ by
\[\mathrm{ev}_M(\sigma) = \sigma(M), \qquad \sigma \in \mathrm{Dens}(W).  \]
The $j$th projection function $p_j(K,\cdot)$ of a convex body $K \in \mathcal{K}(V)$ is now no longer a function on $\mathrm{Gr}_j(V)$ but the section of the
line bundle
\[X_{n-j} = \{(E,l): E \in \mathrm{Gr}_{n-j}(V),\ l \in \mathrm{Dens}(V/E)^*\},  \]
given by
\[p_j(K,E) = \mathrm{ev}_{\mathrm{pr}_E(K)}.  \]
Here, $\mathrm{pr}_E: V \rightarrow V/E$ denotes the natural projection. Similarly, the $(n-j)$th projection function  $p_{n-j}(L,\cdot)$ of $L \in \mathcal{K}(V^*)$
is a section of the line bundle
\[X_j^* = \{(F,l): F \in \mathrm{Gr}_j(V^*),\ l \in \mathrm{Dens}(V^*/F)^*\}.  \]
Let us denote by $C(\mathrm{Gr}_{n-j}(V),X_{n-j})$ and $C(\mathrm{Gr}_j(V^*),X_j^*)$ the spaces of all (continuous) sections of the line bundles $X_{n-j}$ and $X_j^*$, respectively.
Note that the group $\mathrm{SL}(V)$ acts on these vector spaces naturally by left translation. Moreover, using
the annihilator map $\bot: \mathrm{Gr}_{n-j}(V) \rightarrow \mathrm{Gr}_j(V^*)$, it is not difficult to show that the canonical isomorphism
\[\mathrm{Dens}(V/E)^* \cong \mathrm{Dens}(V^*/E^{\bot})^*  \]
induces a canonical $\mathrm{SL}(V)$ equivariant isomorphism between the spaces of sections $C(\mathrm{Gr}_{n-j}(V),X_{n-j})$ and $C(\mathrm{Gr}_j(V^*),X_j^*)$.

For origin-symmetric bodies $K \in \mathcal{K}(V)$ and $L \in \mathcal{K}(V^*)$ with non-empty interior, we may therefore call $K \in \mathcal{K}(V)$ the $j$-projection body of $L$ if
\[p_j(K,\cdot) \cong p_{n-j}(L,\cdot)  \]
with respect to the isomorphism described above. Clearly, this definition coincides with the one given in the introduction if we choose a Euclidean structure
on $V$ and identify $V^*$ with $V$. Furthermore, this invariant formulation immediately implies that the class of $j$-projection bodies is $\mathrm{GL}(V)$ invariant.

\pagebreak

We turn now to classical and new examples of $j$-projection bodies. In the following list, examples $\mathrm{(b)}-\mathrm{(d)}$ were previously considered in the literature.
Examples (e) and (f) are new and based on our main result, Theorem \ref{mainthm1}.

\begin{expls} \label{exp21} \end{expls}

\vspace{-0.3cm}

\begin{enumerate}
\item[(a)] Since, for the Euclidean unit ball $B$ in $\mathbb{R}^n$, we have $\mathrm{vol}_j(B|E^{\bot}) = \kappa_j$ and $\mathrm{vol}_{n-j}(B|E) = \kappa_{n-j}$
for every $1 \leq j \leq n - 1$ and $E \in \mathrm{Gr}_{n-j,n}$, it follows that $B$ is the $j$-projection body of
\[\left ( \frac{\kappa_j}{\kappa_{n-j}} \right )^{\frac{1}{n-j}} B.  \]
Thus, Theorem \ref{thmglninv} and the fact that $(AB)^* = A^{-\mathrm{T}}B$ for every  $A \in \mathrm{GL}(n)$ imply that if $K$ is an origin-symmetric ellipsoid with non-empty interior in $\mathbb{R}^n$, then $K$ is the $j$-projection body of the ellipsoid
\[\left (\frac{\kappa_jV_n(K)}{\kappa_{n-j}\kappa_n} \right )^{\frac{1}{n-j}} K^*.  \]

\item[(b)] McMullen \textbf{\cite{McMullen84}} first proved that the unit cube in $\mathbb{R}^n$,
\[W = \left \{x \in \mathbb{R}^n: -\mbox{$\frac{1}{2}$} \leq x_i \leq \mbox{$\frac{1}{2}$}, i = 1, \ldots, n \right \},  \]
is the $j$-projection body of itself for every $1 \leq j \leq n - 1$. Based on this curious property of $W$, Schnell \textbf{\cite{Schnell94}} deduced the following special case
of Theorem \ref{thmglninv}: If $K$ is an origin-symmetric parallelotope with non-empty interior in $\mathbb{R}^n$, say $K = AW$ with $A \in \mathrm{GL}(n)$,
then $K$ is the $j$-projection body of the parallelotope
\[V_n(K)^{\frac{1}{n-j}} A^{-\mathrm{T}}W.  \]
In particular, if $V_n(K) = 1$, then $A^{-\mathrm{T}}W$ is the $j$-projection body of $K$ for every $1 \leq j \leq n - 1$.
\item[(c)] Let $K \in \mathcal{K}^n_s$ have non-empty interior. Then, by (\ref{suppvol1}) and the definition of the projection body operator $\Pi_{n-1}$ given in
Example \ref{exp17} (a), we have
\[\mathrm{vol}_{n-1}(K|u^{\bot})=h(\Pi_{n-1}K,u) = \mathrm{vol}_{1}\left (\left . \frac{1}{2}\Pi_{n-1}K\right |\mathrm{span}\{u\}\right) \]
for every $u \in \mathbb{S}^{n-1}$,  that is, $K$ is the $(n - 1)$-projection body of $\frac{1}{2}\Pi_{n-1}K$ or, equivalently, $\frac{1}{2}\Pi_{n-1}K$ is the $1$-projection body
of $K$. In particular, the class of $(n-1)$-projection bodies coincides with $\mathcal{K}^n_s$ and, by Minkowski's existence theorem and Cauchy's projection formula, the class of $1$-projection
bodies coincides with the class $\mathcal{Z}_s^n$ of origin-symmetric zonoids in $\mathbb{R}^n$.

\item[(d)] Let $K, L \in \mathcal{K}^n_s$ have non-empty interior. Generalizing a notion introduced by McMullen \textbf{\cite{McMullen87}}, Schneider \textbf{\cite{Schneider96}}
calls $(K,L)$ a \emph{(VP)-pair} if $K$ is the \linebreak $j$-projection body of $L$ for \emph{every} $1 \leq j \leq n - 1$. By the results of McMullen and Schnell described in (b), any pair of
parallelotopes $(AW,A^{-\mathrm{T}}W)$, where $A \in \mathrm{GL}(n)$ and $|\!\det A| = 1$, is an example of a (VP)-pair.

Schneider proved in \textbf{\cite{Schneider96}} that if $K$ and $L$ are polytopes, then they are a (VP)-pair if and only if $K$ is the $j$-projection body of $L$ for $j = 1$ and $j = n - 1$. Moreover, from a result of Weil \textbf{\cite{Weil71}}, Schneider deduced that this holds if and only if $K$ is a direct sum of centrally symmetric polygons and
segments, $V_n(K) = 1$, and $L = \frac{1}{2}\Pi_{n-1}K$.

\item[(e)] Let $K$ be a convex body and $D$ a star body in $\mathbb{R}^{2n}$ and let both be origin-symmetric. Motivated by the identification of $\mathbb{R}^{2n}$ with $\mathbb{C}^n$,
the bodies $K$ and $D$ are called \emph{complex} if they are invariant with respect to any coordinate-wise two-dimensional rotation (see, e.g., \textbf{\cite{KolPaoZym13}} for details).
Note that origin-symmetric complex convex bodies in $\mathbb{R}^{2n}$ correspond precisely to the unit balls of complex norms on $\mathbb{C}^n$.

For a unit vector $u \in \mathbb{C}^n$, let $H_u = \left \{z \in \mathbb{C}^n: \langle u,z \rangle = \sum_{k=1}^n u_k\overline{z}_k = 0 \right \}$
denote the complex hyperplane perpendicular to $u$. Under the standard mapping from $\mathbb{C}^n$ to $\mathbb{R}^{2n}$, the hyperplane $H_u$ becomes a
$2n-2$ dimensional subspace of $\mathbb{R}^{2n}$ which is orthogonal to the vectors
\[u=(u_{11},u_{12},\ldots,u_{n1},u_{n2}) \qquad \mbox{and} \qquad u_* = (-u_{12},u_{11},\ldots,-u_{n2},u_{n1}).  \]

\noindent {\bf Definition}\! (\!\!\!\textbf{\cite{KolPaoZym13}}) \emph{Let $D$ and $M$ be origin-symmetric complex star bodies in $\mathbb{R}^{2n}$.
Then $D$ is called the \emph{complex intersection body} of $M$ if
\[\mathrm{vol}_{2}(D \cap H_u^{\bot}) = \mathrm{vol}_{2n-2}(M \cap H_u)  \]
for every $u \in \mathbb{S}^{2n-1}$. The class of  complex intersection bodies is the closure in the radial metric of all complex intersection bodies of star bodies.}

Koldobsky, Paouris, and Zymonopoulou \textbf{\cite{KolPaoZym13}} proved that the class of complex intersection bodies coincides with the class of $2$-intersection bodies which are complex.
Moreover, they showed that complex intersection bodies of convex bodies are also convex. Motivated by these results, we define complex projection bodies as follows
(see \textbf{\cite{AbardiaBernig11}}, for a different notion of complex projection bodies).

\noindent {\bf Definition} \emph{Let $K$ and $L$ be origin-symmetric complex convex bodies in $\mathbb{R}^{2n}$.
Then $K$ is called the \emph{complex projection body} of $L$ if
\[\mathrm{vol}_{2}(K|H_u^{\bot}) = \mathrm{vol}_{2n-2}(L|H_u)  \]
for every $u \in \mathbb{S}^{2n-1}$. The class of  complex projection bodies is the closure in the Hausdorff metric of all complex projection bodies of convex bodies.}

\pagebreak

Clearly, if $K \in \mathcal{K}^{2n}_s$ is a complex $2$-projection body of $L$, then $K$ is a the complex projection body of $L$.
We do not know if the converse also holds. However, if $D \in \mathcal{K}^{2n}_s$ is complex and of class $C^{\infty}_+$, then, by Proposition~\ref{guanma}, the function $\rho(D,\cdot)^2 \in C_e^{\infty}(\mathbb{S}^{n-1})$ is the density of the area measure of order $2$ of a complex convex body $K \in \mathcal{K}^{2n}_s$. Consequently,
by Theorems \ref{koldmain} and \ref{mainthm1} and Proposition \ref{guanma}, if $D$ is a complex intersection body, then $K$ is a complex $2$-projection body which, in turn,
is a complex projection body.

\item[(f)] Finally, we consider strictly convex bodies of revolution $K_{\lambda} \in \mathcal{K}^n_s$ whose area measures of order $1 \leq j \leq n - 2$ have a density of the form
\begin{equation} \label{fireyexp}
s_j(K_{\lambda},\bar{e} \cdot .\,) = 1 + \lambda P_2^n(\bar{e} \cdot.\,) = 1 + \frac{\lambda}{n-1}(n(\bar{e}\cdot .\,)^2 - 1),
\end{equation}
where $P_2^n$ denotes the Legendre polynomial of dimension $n$ and degree~$2$.
In order to determine all admissible $\lambda$ in (\ref{fireyexp}), we use Theorem \ref{firey_area_measures}.
Clearly, condition (i) of Theorem \ref{firey_area_measures} is satisfied for all $\lambda \in \mathbb{R}$. However, since
\[\int_t^1 \xi\, s_j(K_{\lambda},\xi)(1-\xi^2)^{\frac{n-3}{2}}d\xi = \frac{(1-t^2)^{\frac{n-1}{2}}}{n^2-1}(\lambda(nt^2+1)+n+1),  \]
it is not difficult to show that conditions (ii) and (iii) of Theorem \ref{firey_area_measures} are satisfied if and only if
\begin{equation} \label{lambdarange1}
\lambda \in \left (-1,\frac{j(n+1)}{2n-j} \right ).
\end{equation}
Now, we want to determine which of the bodies $K_{\lambda}$ are $j$-projection bodies.
To this end, note that, by (\ref{multsphfour}), we have
\[\mathbf{F}_{-j}s_j(K_{\lambda},\bar{e} \cdot .\,) = \frac{\pi^{\frac{n}{2}}2^{n-j}\Gamma\left (\frac{n-j}{2} \right )}{\Gamma\left (\frac{j}{2}\right )}
\left (1- \lambda \frac{n-j}{j} P_2^n(\bar{e} \cdot.\,)  \right ).  \]
Hence, by Theorem \ref{mainthm1} and (\ref{lambdarange1}), $K_{\lambda}$ is a $j$-projection body if and only if
\[\lambda \in \left (-1,\frac{j}{n-j} \right ).  \]
This shows, in particular, that for $j < n - 1$ the class of $j$-projection bodies is a proper subset of $\mathcal{K}^n_s$.
\end{enumerate}

In the final part of this section, we want to prove two more basic properties of $j$-projection bodies. The first one is a generalization
of the well known fact that Minkowski's projection body operator $\Pi_{n-1}$ maps polytopes to polytopes. Note that, by Example \ref{exp21} (c), this implies
that $1$- and $(n - 1)$-projection bodies of polytopes are polytopes. As part of the following result we extend this observation to all $j \in \{1, \ldots, n - 1\}$.

\begin{theorem} Let $1 \leq j \leq n - 1$ and let $P$ and $Q$ be origin-symmetric convex bodies with non-empty interior in $\mathbb{R}^n$.
If $P$ is the $j$-projection body of $Q$, then
\begin{equation} \label{maineqpoly}
S_j(P,\mathbb{S}^{n-1} \cap E) = \frac{\kappa_{n-j}}{\kappa_j}\,S_{n-j}(Q,\mathbb{S}^{n-1} \cap E^{\bot})
\end{equation}
for every $E \in \mathrm{Gr}_{n-j,n}$. Moreover, if $P$ is a polytope, then so is $Q$ and
$P$ has a $j$-face parallel to $E$ if and only if $Q$ has an $(n-j)$-face parallel to $E^{\bot}$.
\end{theorem}

\noindent {\it Proof}.  In order to prove (\ref{maineqpoly}), let $E \in \mathrm{Gr}_{n-j,n}$ be an arbitrary but fixed subspace.
For $\varepsilon > 0$, let $f_{\varepsilon} \in C([0,1])$ be monotone increasing with $\mathrm{supp}\,f_{\varepsilon} \subseteq [1-2\varepsilon,1]$
and such that $f_{\varepsilon} \equiv 1$ on $[1-\varepsilon,1]$ and define $g_{\varepsilon}^E \in C(\mathrm{Gr}_{j,n})$, by
\[g^{E}_{\varepsilon}(F) = f_{\varepsilon}(|\cos(E,F)|). \]
Note that
\[(\mathrm{R}_{j,1}g^{E}_{\varepsilon})(u) = \int_{\mathrm{Gr}_{j,n}^u} f_{\varepsilon}(|\cos(E,F)|)\,d\nu_j^u(F)   \]
depends only on $|\cos(E,u)|$. In particular, $\mathrm{R}_{j,1}g^{E}_{\varepsilon}$ is constant on $\mathbb{S}^{n-1} \cap E$.
Consequently, we can replace $g^{E}_{\varepsilon}$, if necessary, by a positive multiple such that $(\mathrm{R}_{j,1}g^{E}_{\varepsilon})(u)=1$
whenever $u \in \mathbb{S}^{n-1} \cap E$. Next, we want to show that
\begin{equation} \label{radongeps}
 \lim_{\varepsilon \to 0} (\mathrm{R}_{j,1}g^{E}_{\varepsilon})(u) = \left \{ \begin{array}{ll} 1 & \mbox{for } u \in \mathbb{S}^{n-1} \cap E, \\
 0 & \mbox{for } u \notin \mathbb{S}^{n-1} \cap E. \end{array} \right.
\end{equation}
To this end, observe that $|\cos(E,F)| \leq |\cos (E,u)|$ whenever $u \in E$. Thus, by the monotonicity of $f_{\varepsilon}$, we have
\[ (\mathrm{R}_{j,1}g^{E}_{\varepsilon})(u) \leq  f_{\varepsilon}(|\cos (E,u)|).\]
Hence, by the definition of $f_{\varepsilon}$, for every $u \notin \mathbb{S}^{n-1} \cap E$ there exists $\varepsilon_{u}>0$ such that
$(\mathrm{R}_{j,1}g^{E}_{\varepsilon})(u)=0$ for every $\varepsilon \leq \varepsilon_{u}$ which completes the proof of (\ref{radongeps}).

The same arguments used to prove (\ref{radongeps}) together with the fact that $\mathrm{SO}(n)$ acts transitively on $\mathrm{Gr}_{n-j,n}$, show that there exists a positive constant $c \in \mathbb{R}$,
independent of $E$, such that
\begin{equation} \label{radongepsbot}
 \lim_{\varepsilon \to 0} (\mathrm{R}_{n-j,1}(g^{E}_{\varepsilon})^{\bot})(u) = \left \{ \begin{array}{ll} c & \mbox{for } u \in \mathbb{S}^{n-1} \cap E^{\bot}, \\
 0 & \mbox{for } u \notin \mathbb{S}^{n-1} \cap E^{\bot}. \end{array} \right.
\end{equation}

Now, since $P$ is the $j$-projection body of $Q$ and $\mathrm{R}_{j,i}$ is the adjoint of $\mathrm{R}_{i,j}$, it follows from
Theorem \ref{mainthm1} and Proposition \ref{spherfourradon} that
\begin{eqnarray*}
\int_{\mathbb{S}^{n-1}}(\mathrm{R}_{j,1}g^{E}_{\varepsilon})(u)\,dS_j(P,u) & = & \int_{\mathrm{Gr}_{j,n}}g^{E}_{\varepsilon}(F)\,d(\mathrm{R}_{1,j}S_j(P,\cdot))(F) \\
& = & \frac{\kappa_{n-j}}{\kappa_j} \int_{\mathrm{Gr}_{n-j,n}}g^{E}_{\varepsilon}(F^{\bot})\,d(\mathrm{R}_{1,n-j}S_{n-j}(Q,\cdot))(F) \\
& = & \frac{\kappa_{n-j}}{\kappa_j} \int_{\mathbb{S}^{n-1}}(\mathrm{R}_{n-j,1}(g^{E}_{\varepsilon})^{\bot})(u)\,dS_{n-j}(Q,u).
\end{eqnarray*}
Letting $\varepsilon \to 0$, (\ref{radongeps}), (\ref{radongepsbot}), and the dominated convergence theorem yield
\begin{equation*}
S_j(P,\mathbb{S}^{n-1} \cap E) = \frac{c\,\kappa_{n-j}}{\kappa_j}\, S_{n-j}(Q,\mathbb{S}^{n-1} \cap E^{\bot}),
\end{equation*}
where the positive constant $c$ is the same as in (\ref{radongepsbot}) and does not depend on $P$ and $Q$.
Using Example \ref{exp21} (b) and taking $P = W = Q$, shows, by (\ref{sjpoly}), that $c = 1$ which completes the proof of (\ref{maineqpoly}).

Now assume that $P$ is a polytope and recall that, by (\ref{sjpoly}), $S_j(P,\cdot)$ is concentrated on the union of finitely many $n-1-j$ dimensional great spheres.
Let $\mathcal{G}^{n-j}(P)$ denote the \emph{finite} set of subspaces $E \in \mathrm{Gr}_{n-j,n}$ such that $S_j(P,\mathbb{S}^{n-1} \cap E) > 0$.
Summing (\ref{maineqpoly}) over all $E \in \mathcal{G}^{n-j}(P)$ yields on one hand
\[\frac{\kappa_{n-j}}{\kappa_j}\sum_{E\in \mathcal{G}^{n-j}(P)}\!\!\! S_{n-j}(Q,\mathbb{S}^{n-1} \cap E^{\bot}) =
\sum_{E\in \mathcal{G}^{n-j}(P)}\!\!\! S_j(P,\mathbb{S}^{n-1} \cap E) = S_j(P,\mathbb{S}^{n-1}). \]
On the other hand, since $Q$ is the $(n-j)$-projection body of $P$, Theorem \ref{mainthm1} and (\ref{multsphfour}) imply that
\[S_j(P,\mathbb{S}^{n-1}) = \frac{ja_0^n[\mathbf{F}_{j-n}]}{(2\pi)^j(n-j)} S_{n-j}(Q,\mathbb{S}^{n-1}) = \frac{\kappa_{n-j}}{\kappa_j} S_{n-j}(Q,\mathbb{S}^{n-1}).  \]
Consequently, since $S_{n-j}(Q,\cdot)$ vanishes on great spheres of dimension $d < j-1$, we have
\[S_{n-j}(Q,\mathbb{S}^{n-1}) = \sum_{E\in \mathcal{G}^{n-j}(P)}\!\!\! S_{n-j}(Q,\mathbb{S}^{n-1} \cap E^{\bot}) = S_{n-j} \left(Q,\!\!\! \bigcup_{E \in \mathcal{G}^{n-j}(P)}\!\!\! (\mathbb{S}^{n-1} \cap E^{\bot}) \right).  \]
This shows that $S_{n-j}(Q, \cdot)$ is concentrated on a finite union of $j-1$ dimensional great spheres.
An application of Proposition \ref{sjpolytope} finishes the proof.
\hfill $\blacksquare$

\vspace{0.3cm}

It is well known that for every pair of convex bodies $K, L \in \mathcal{K}^n$,
\begin{equation} \label{pimixedid}
V(K[n-1],\Pi_{n-1}L) =  V(L[n-1],\Pi_{n-1}K).
\end{equation}
This basic mixed volume identity for Minkowski's projection body operator $\Pi_{n-1}$ and its variants for other Minkowski valuations have found
numerous applications (see, e.g., \textbf{\cite{AbardiaBernig11, ABS11, GooWei12, GooZha98, Schuster10}}).
In view of Example \ref{exp21} (c), our final result of this section provides a generalization of (\ref{pimixedid}) in the context of $j$-projection bodies.

\begin{theorem} Let $1 \leq j \leq n - 1$ and let $K_i$, $L_i$, $i = 1, 2$, be origin-symmetric convex bodies with non-empty interior in $\mathbb{R}^n$.
If $K_1$ is the $j$-projection body of $L_1$ and $K_2$ is the $(n-j)$-projection body of $L_2$, then
\[V(K_1[j],K_2[n-j]) = V(L_1[n-j],L_2[j]).  \]
\end{theorem}

\noindent {\it Proof}. Consider the valuations $\phi \in \mathbf{Val}_{n-j}^+$ and $\psi \in \mathbf{Val}_j^+$, defined by
\[\phi(K) = V(K_1[j], K[n-j]) \qquad \mbox{and} \qquad  \psi(K) = V(L_1[n-j],K[j]).\]
Then, by a well known relation between projection functions and mixed volumes (see, e.g., \textbf{\cite[\textnormal{Theorem 5.3.1}]{Schneider14}}),
the Klain functions $\mathrm{Kl}_{n-j}\phi \in C(\mathrm{Gr}_{n-j,n})$ and $\mathrm{Kl}_j\psi \in C(\mathrm{Gr}_{j,n})$ are given by
\[\mathrm{Kl}_{n-j}\phi(E) = \binom{n}{j}^{-1} \mathrm{vol}_j(K_1|E^{\bot}) \quad \,\, \mbox{and} \quad \,\,  \mathrm{Kl}_j\psi(F) = \binom{n}{j}^{-1}\!\! \mathrm{vol}_{n-j}(L_1|F^{\bot}).\]
Therefore, since $K_1$ is the $j$-projection body of $L_1$, we have
\begin{equation} \label{mixedvolid1}
\mathrm{Kl}_{n-j}\phi = (\mathrm{Kl}_j\psi)^{\bot}.
\end{equation}

Now, assume that $\phi$ and $\psi$ are smooth. Then, by (\ref{mixedvolid1}) and (\ref{deffourier}), $\psi = \mathbb{F}\phi$. Moreover,
(as already explained after the proof of Theorem \ref{mainequiv}) it follows from Theorem \ref{klaincoscrof} and (\ref{klaincrofton}) that
\begin{equation} \label{mixedvolid2}
\phi(K_2) = V(K_1[j], K_2[n-j]) = V(L_1[n-j],L_2[j]) = \psi(L_2)
\end{equation}
which is the desired relation.

If $\phi$ and $\psi$ are not smooth, but merely continuous, then a recent extension of Alesker and Faifman \textbf{\cite[\textnormal{Propositions 4.4 and 4.5}]{AleskerFaifman}} of the Klain and Crofton maps as well as Theorem \ref{klaincoscrof} and (\ref{klaincrofton}) to generalized valuations (which include, in particular, continuous valuations) implies that (\ref{mixedvolid2}) still follows from
(\ref{mixedvolid1}). \hfill $\blacksquare$

\vspace{1cm}

\centerline{\large{\bf{ \setcounter{abschnitt}{6}
\arabic{abschnitt}. From $j$-intersection bodies to $j$-projection bodies}}}

\reseteqn \alpheqn \setcounter{theorem}{0}

\vspace{0.6cm}

In this final section, we first recall the definition of the class of convex bodies $\mathcal{K}_s^n(j)$ and their dual analogs, the class of $j$-Busemann-Petty star bodies (also called generalized $j$-intersection bodies). Then, we relate these two classes as well as the classes of $j$-intersection bodies and $j$-projection bodies via a generalization of the duality transform introduced at the end of Example~\ref{exp21}~(e). Finally, we prove a dual analog of a recent result of Milman on the relation between the classes of $j$-intersection bodies and $j$-Busemann-Petty bodies
and state two open problems on the class of $j$-projection bodies.

For $1 \leq j \leq n - 1$, let $\mathcal{P}_s^n(j)$ denote the class of all (origin-symmetric) \linebreak $j$-projection bodies in $\mathbb{R}^n$ and recall that $\mathcal{K}_s^n(j)$ is the class
of all origin-symmetric convex bodies $K \in \mathcal{K}^n$ such that
\begin{equation*}
\mathrm{vol}_j(K|\,\cdot\,) = \mathrm{C}_j\varrho_j(K,\cdot)
\end{equation*}
for some \emph{non-negative} Borel measure $\varrho_j(K,\cdot)$ on $\mathrm{Gr}_{j,n}$.

\pagebreak

Using the Cauchy-Kubota formula (\ref{radsurfcos}) as well as (\ref{cifbot}) and (\ref{rijbot}), it follows that $K \in \mathcal{K}_s^n(j)$ if and only if
\[\mathrm{C}_1S_j(K,\cdot) = \frac{2\kappa_{n-1}}{\kappa_j}(\mathrm{R}_{n-j,1} \circ \mathrm{C}_{n-j})\varrho_j^{\bot}(K,\cdot)\]
which, by the composition rule for Radon and cosine transforms (see, e.g., \textbf{\cite{GooZha98}}), is equivalent to
\begin{equation} \label{projgensurfarea}
\mathrm{C}_1S_j(K,\cdot) = {n \choose j}^{-1}\!\!n\kappa_{n-j} (\mathrm{C}_1 \circ \mathrm{R}_{n-j,1})\varrho_j^{\bot}(K,\cdot).
\end{equation}
Consequently, the injectivity of the spherical cosine transform $\mathrm{C}_1$ implies that the class $\mathcal{K}_s^n(j)$ consists precisely of those origin-symmetric
$K \in \mathcal{K}^n$ for which
\begin{equation} \label{kjsurfarea}
S_j(K,\cdot) = \mathrm{R}_{n-j,1} \mu_j(K,\cdot)
\end{equation}
for some \emph{non-negative} Borel measure $\mu_j(K,\cdot)$ on $\mathrm{Gr}_{n-j,n}$.

The dual analog of the class $\mathcal{K}_s^n(j)$ was introduced by Zhang in 1996.

\vspace{0.3cm}

\noindent {\bf Definition}  (\!\!\textbf{\cite{zhang96}}) \emph{Suppose that $1 \leq j \leq n - 1$. An origin-symmetric star body $D$ in $\mathbb{R}^n$ is called
a $j$-Busemann-Petty body if
\begin{equation} \label{defbusepetty}
\rho(D,\cdot)^j = \mathrm{R}_{n-j,1} \nu_j(D,\cdot)
\end{equation}
for some non-negative Borel measure $\nu_j(D,\cdot)$ on $\mathrm{Gr}_{n-j,n}$.}

\vspace{0.3cm}

For $1 \leq j \leq n - 1$, let $\mathcal{I}_s^n(j)$ denote the class of all (origin-symmetric) \linebreak $j$-intersection bodies in $\mathbb{R}^n$ and
let $\mathcal{BP}_s^n(j)$ denote the class of (origin-symmetric) \linebreak $j$-Busemann-Petty bodies in $\mathbb{R}^n$. From their definition and Theorem \ref{koldmain} it follows easily
that
\begin{equation} \label{inclusionint}
\mathcal{I}_s^n(1) = \mathcal{BP}_s^n(1) \qquad \mbox{and} \qquad \mathcal{I}_s^n(n - 1) = \mathcal{BP}_s^n(n-1) = \mathcal{S}^n_s,
\end{equation}
where here and in the following $\mathcal{S}_s^n$ denotes the class of origin-symmetric star bodies in $\mathbb{R}^n$.
Recall that $\mathcal{I}_s^n(1)$ coincides with Lutwak's intersection bodies.

The discovery and importance of the class $\mathcal{BP}_s^n(j)$ is due to their connection to the $j$-codimensional Busemann-Petty problem which asks
whether the volume of a \emph{convex} body $K \in \mathcal{K}_s^n$ is smaller than that of another body $L \in \mathcal{K}^n_s$ if
\[\mathrm{vol}_{n-j}(K \cap E) \leq \mathrm{vol}_{n-j}(L \cap E)  \]
for all $E \in \mathrm{Gr}_{n-j,n}$. Zhang \textbf{\cite{zhang96}} showed that a positive answer to this problem is equivalent to whether all origin-symmetric convex bodies in $\mathbb{R}^n$ belong to the class
$\mathcal{BP}_s^n(j)$. Subsequently, Bourgain and Zhang \textbf{\cite{BourgainZhang}} proved that the answer is negative for $j < n - 3$ but the cases $j = n - 3$ and $j = n - 2$ remained open.
This was later reproved by Koldobsky \textbf{\cite{Koldobsky00}} who also first considered the relationship between the two types of generalizations of Lutwak's intersection bodies,
$\mathcal{I}_s^n(j)$ and $\mathcal{BP}_s^n(j)$, and proved that for all $1 \leq j \leq n - 1$,
\begin{equation} \label{bpsubseti}
\mathcal{I}_s^n(1) \subseteq \mathcal{BP}_s^n(j) \subseteq \mathcal{I}_s^n(j).
\end{equation}

\pagebreak

Koldobsky also asked whether, in fact, $\mathcal{BP}_s^n(j) = \mathcal{I}_s^n(j)$ holds not just for $j = 1$ and $j = n - 1$ but for all $2 \leq j \leq n - 2$ as well.
If this were true, a positive answer to the $j$-codimensional Busemann-Petty problem for $j \geq n - 3$ would follow, since $\mathcal{K}_s^n \subseteq \mathcal{I}_s^n(j)$ for those values of $j$.
However, Milman gave the following negative answer to Koldobsky's question.

\begin{theorem} \label{milmanmain} \emph{(\!\!\textbf{\cite{Milman08}})}  Suppose that $n \geq 4$ and $2 \leq j \leq n - 2$. Then there exists a smooth
star body of revolution $D \in \mathcal{I}_s^n(j)$ such that $D \not \in \mathcal{BP}_s^n(j)$.
\end{theorem}

Note that Theorem \ref{milmanmain} did not resolve the open cases of the $j$-codimensional Busemann-Petty problem since the body $D$ is not necessarily convex.

\vspace{0.2cm}

Motivated by the formal analogy between Theorems \ref{koldmain} and \ref{mainthm1} as well as the definitions (\ref{kjsurfarea}) and (\ref{defbusepetty}), we define now two 'duality' transforms on smooth convex bodies using
Proposition \ref{guanma} of Guan and Ma, thereby extending the map which already appeared at the end of Example \ref{exp21} (e).
To this end, let $\mathcal{K}^{n,\infty}_s$ denote the subset of $\mathcal{K}^n_s$ consisting of convex bodies of class $C^{\infty}_+$.

\vspace{0.3cm}

\noindent {\bf Definition}  \emph{For $1 \leq j \leq n - 1$, the map $\mathrm{P}_j: \mathcal{K}^{n,\infty}_s \rightarrow \mathcal{K}^{n,\infty}_s$ is defined by
\[s_j(\mathrm{P}_jK,\cdot) = \rho(K,\cdot)^j.  \]
The map $\mathrm{I}_j: \mathcal{K}^{n,\infty}_s \rightarrow \mathcal{S}^n_s$ is defined by}
\[\rho(\mathrm{I}_jK,\cdot)^j = s_j(K,\cdot).  \]

\vspace{0.2cm}

Clearly, the map $\mathrm{I}_j$ is a left inverse of $\mathrm{P}_j$, that is $\mathrm{I}_j \circ \mathrm{P}_j = \mathrm{id}$.
Moreover, the following immediate consequences of  Theorems \ref{koldmain} and \ref{mainthm1} and definitions (\ref{kjsurfarea}) and (\ref{defbusepetty}),
show that these maps are closely related to the various notions of intersection and projection bodies.

\begin{coro} Suppose that $1 \leq j \leq n - 1$ and let $K, L \in \mathcal{K}^{n,\infty}_s$.
\begin{enumerate}
\item[(a)] If $K$ is $j$-intersection body of $L$, then $\mathrm{P}_jK$ is $j$-projection body of $\mathrm{P}_{n-j}L$.
\item[(b)] If $K$ is $j$-projection body of $L$, then $\mathrm{I}_jK$ is $j$-intersection body of $\mathrm{I}_{n-j}L$.
\item[(c)] If $K \in \mathcal{BP}_s^n(j)$, then $\mathrm{P}_jK \in \mathcal{K}_s^n(j)$.
\item[(d)] If $K \in \mathcal{K}_s^n(j)$, then $\mathrm{I}_jK \in \mathcal{BP}_s^n(j)$.
\end{enumerate}
\end{coro}

In view of the duality relations and various analogies between results on $j$-intersection bodies and $j$-projection bodies we have encountered so far, it is natural to ask whether there is a relation similar to (\ref{bpsubseti}) between the classes $\mathcal{K}_s^n(j)$ and $\mathcal{P}_s^n(j)$. Recall from Section 3 and Example \ref{exp21} (c) that
\begin{equation*}
\mathcal{K}_s^n(1) = \mathcal{P}_s^n(1) = \mathcal{Z}_s^n \qquad \mbox{and} \qquad \mathcal{K}_s^n(n - 1) = \mathcal{P}_s^n(n-1) = \mathcal{K}^n_s,
\end{equation*}
which is the dual analog of (\ref{inclusionint}), and $\mathcal{Z}_s^n \subseteq \mathcal{K}_s^n(j)$ for all $1 \leq j \leq n - 1$.

\pagebreak

\noindent {\bf Open Problems}  \emph{Suppose that $2 \leq j \leq n - 2$.}
\begin{enumerate}
\item[(a)] Is it true that $\mathcal{K}_s^n(j) \subseteq \mathcal{P}_s^n(j)$?
\item[(b)] Is it true that $\mathcal{Z}_s^n(j) \subseteq \mathcal{P}_s^n(j)$?
\end{enumerate}

For the one-parameter family of $j$-projection bodies $K_{\lambda}$ constructed in Example \ref{exp21} (f) it is possible to show that
both problems have a positive answer. In fact, the entire family $K_{\lambda}$ is also contained in $\mathcal{K}_s^n(j)$. However, with our final result
we show that this is not always the case by proving a partial analog of Theorem \ref{milmanmain}.

\begin{theorem} There exists a smooth convex body of revolution $K \in \mathcal{P}_s^4(2)$ such that $K \not \in \mathcal{K}_s^4(2)$.
\end{theorem}

\noindent {\it Proof}. For $\varepsilon > 0$ and $t \in [-1,1]$, let
\[ s_{\varepsilon}(t) = 1 + \varepsilon +  \frac{5}{2}P^4_4( t) = \frac{3}{2} + \varepsilon - 6t^2 + 8 t^4,\]
where $P_4^4$ denotes the Legendre polynomial of dimension $4$ and degree $4$. We first want to prove
that for every $\varepsilon > 0$, there exists a strictly convex body of revolution $K_{\varepsilon} \in \mathcal{K}_s^4$ whose
area measure of order $2$ has a density of the form
\[s_2(K_{\varepsilon},\bar{e}\cdot.\,) = s_{\varepsilon}(\bar{e}\cdot.\,).  \]
To this end, we will show that $s_{\varepsilon}(\bar{e}\cdot . \,)$ satisfies conditions (i), (ii), and (iii) of Theorem \ref{firey_area_measures} with $j=2$ and $n=4$.
Note that since $s_{\varepsilon}$ is an \emph{even} polynomial which is strictly positive for every $\varepsilon > 0$, conditions (i) and (ii) hold trivially.
However, from
\[\int_t^1 \xi\,s_{\varepsilon}(\xi)(1-\xi^2)^{\frac{1}{2}}\,d\xi = \frac{(1-t^2)^{\frac{3}{2}}}{42}(13 + 14\varepsilon - 12t^2 + 48t^4)  \]
it follows by a straightforward calculation that also condition (iii) is satisfied.

Next observe that, by (\ref{multsphfour}), $a_{2k}^4[\mathbf{F}_{-2}] = (2\pi)^2(-1)^k$. Thus, by Theorem \ref{mainthm1}, $K_{\varepsilon}$ is a $2$-projection body
of a dilate of itself. It remains to show that $K_{\varepsilon} \notin \mathcal{K}_s^4(2)$ for sufficiently small $\varepsilon$. To this end, note that
there exists a \emph{unique spherical} function $\varrho_2(K_{\varepsilon},\cdot) \in C^{\infty}(\mathrm{Gr}_{j,n})^{\mathrm{sph}}$ such that
\begin{equation} \label{lastproof1}
\mathrm{vol}_2(K_{\varepsilon}|\,\cdot\,) = \mathrm{C}_2 \varrho_2(K_{\varepsilon},\cdot).
\end{equation}
Indeed, by (\ref{projgensurfarea}), the function $\varrho_2(K_{\varepsilon},\cdot)$ is given by
\[\varrho_2(K_{\epsilon}, \cdot) = \frac{3}{2\pi} (\bot_{\ast} \circ \mathrm{R}_{2,3})^{-1} s_{\varepsilon}(\bar{e}\cdot \, . \,).\]
This is well defined since $s_2(K_\varepsilon, \cdot)$ is smooth and $\mathrm{SO}(n-1)$-invariant and therefore, by definition, spherical.
The uniqueness follows from the injectivity of $\mathrm{C}_2$ on spherical functions. We conclude that, since $K_{\varepsilon}$ is $\mathrm{SO}(n-1)$ invariant and $\mathrm{C}_2$ commutes with rotations,
the $\mathrm{SO}(n-1)$ symmetrization of any measure satisfying (\ref{lastproof1}) must coincide with $\varrho_2(K_{\varepsilon},\cdot)$. Hence, in order to prove that
$K_{\varepsilon} \notin \mathcal{K}_s^4(2)$, it suffices to show that $(\bot_{\ast} \circ \mathrm{R}_{2,3})^{-1} s_0(\bar{e}\cdot \, . \,)$ attains negative values.

Using the spherical Radon transform $\mathrm{R}:= \bot_* \circ \mathrm{R}_{1,n-1}$ and the composition rules for Radon transforms, it follows that
\[(\bot_{\ast} \circ \mathrm{R}_{2,3})^{-1} s_0(\bar{e}\cdot \, . \,) =\left(\mathrm{R}_{1,2} \circ \mathrm{R}^{-1} \right) s_0(\bar{e}\cdot \, . \,).\]
Hence, since $a_0^4[\mathrm{R}] = 1$ and $a_4^4[\mathrm{R}] = \frac{1}{5}$  (see, e.g., \textbf{\cite[\textnormal{Lemma 3.4.7}]{Groemer}}), we must show that
\[ (\mathrm{R}_{1,2} \circ \mathrm{R}^{-1} ) \, s_0( \bar{e} \cdot \, . \,) =  \mathrm{R}_{1,2} \left( 1  + \frac{25}{2}P_4^4(\bar{e} \cdot \, . \,) \right) \]
attains negative values. Now, by \textbf{\cite[\textnormal{Corollary 3.3}]{Milman08}}, we have for $f \in C[0,1]$,
\[\left( R_{1,2} f(\bar{e} \cdot . \,) \right) (E) = \frac{2}{\pi} \int_0^1 f(|\cos(E, \bar{e})|\, t) (1-t^2)^{-\frac{1}{2}} \, dt .\]
This means that we have to find a $\xi \in [0,1]$ such that
\[\int_0^1 \left(1  + \frac{25}{2} P_4^4(\xi t) \right) (1-t^2)^{-\frac{1}{2}} \, dt =  \frac{15 \pi}{2} \xi^2(\xi^2 - 1) + \frac{7\pi}{4}  < 0.\]
Clearly, one possible choice is given by $\xi = \frac{1}{\sqrt{2}}$. \hfill $\blacksquare$

\vspace{1cm}

\noindent {{\bf Acknowledgments} The work of the authors was
supported by the European Research Council (ERC), Project number: 306445, and the Austrian Science Fund (FWF), Project number:
Y603-N26.

\noindent Vienna University of Technology \par \noindent Institute of Discrete
Mathematics and Geometry \par \noindent Wiedner Hauptstra\ss e 8--10/1047
\par \noindent A--1040 Vienna, Austria

\vspace{0.2cm}

\par \noindent felix.dorrek@tuwien.ac.at \par \noindent franz.schuster@tuwien.ac.at

\end{document}